\documentclass{cocv}

\usepackage{amssymb}
\usepackage[english]{babel}
\usepackage[utf8]{inputenc}
\usepackage{graphicx}
\usepackage{tikz}
\usepackage{pgfplots}
\usepackage{enumitem}
\usepackage{xfrac}
\usepackage{mathtools}
\usepackage{url}
\pgfplotsset{compat=1.18}

\newcommand{\dint}[1]{\;\mathrm{d} #1 }
\newcommand{\R}{\mathbb{R}}

\newcommand{\N}{\mathbb{N}}

\newcommand{\fLp}[2]{\mathrm{L}^{ #1 }\left( #2 \right)}
\newcommand{\fLpw}[3]{\mathrm{L}_{#2}^{ #1 }\left( #3 \right)}
\newcommand{\Lp}[2]{\mathrm{L}^{ #1 }}
\newcommand{\Lpw}[3]{\mathrm{L}_{#2}^{ #1 }}

\newcommand{\fWp}[2]{\mathrm{W}^{1, #1 }\left( #2 \right)}
\newcommand{\fWpw}[3]{\mathrm{W}^{1,#1}_{ #2 }\left( #3 \right)}
\newcommand{\Wp}[2]{\mathrm{W}^{1, #1 }}
\newcommand{\Wpw}[3]{\mathrm{W}^{1,#1}_{ #2 }}

\newcommand{\fNWpw}[3]{\mathrm{W}^{-1,#1}_{ #2 }\left( #3 \right)}

\newcommand{\NWpw}[3]{\mathrm{W}^{-1,#1}_{ #2 }}

\newcommand{\CcInfty}[1]{\mathcal{C}_c^\infty\left( #1 \right) }

\newcommand{\inner}[3][]{ \left< #2, #3\right>_{ #1 } }
\newcommand{\einner}[1]{ \left< \cdot, \cdot \right>_{ #1 } }

\renewcommand{\det}[1]{\operatorname{det}\big( #1 \big )}

\DeclareMathOperator*{\Range}{range}
\DeclareMathOperator*{\Image}{im\,}

\DeclareMathOperator*{\clSpan}{\overline{span}\;}
\DeclareMathOperator*{\esssup}{ess\,sup}

\DeclareMathOperator*{\supp}{supp}
\DeclareMathOperator*{\wslim}{w^\ast-\operatorname{lim}}

\mathchardef\ordinarycolon\mathcode`\:
\mathcode`\:=\string"8000
\begingroup \catcode`\:=\active
\gdef:{\mathrel{\mathop\ordinarycolon}}
\endgroup

\begin{document}

	\title{Unifying HJB and Riccati equations: A Koopman operator approach to nonlinear optimal control}
	\runningtitle{Unifying HJB and Riccati equations}
		\author{Bernhard Höveler}
	\address{RWTH Aachen University, Institut für Geometrie und Praktische Mathematik}
	\author{Tobias Breiten}
	\address{Technische Universität Berlin, Institute of Mathematics}

	\subjclass{49L12, 47D06, 93C10, 93D20}
	\keywords{nonlinear optimal control, Hamilton--Jacobi--Bellman equation, Koopman operator, operator Riccati equation, low-rank approximation}

	\begin{abstract} 
		This paper proposes an operator-theoretic framework that recasts the minimal value function of a nonlinear optimal control problem as an abstract bilinear form on a suitable function space. The resulting bilinear form is shown to satisfy an operator equation with quadratic nonlinearity obtained by formulating the Lyapunov equation for a Koopman lift of the optimal closed-loop dynamics to an infinite-dimensional state space. It is proved that, under regularity assumptions, the minimal value function admits a rapidly convergent sum-of-squares expansion, a direct consequence of the fast spectral decay of the bilinear form. The framework thereby establishes a natural link between the Hamilton--Jacobi--Bellman equation and a Riccati-like operator equation.
	\end{abstract}

	\maketitle
	
	\section{Introduction}
	\label{sec:intro}
	We consider nonlinear control-affine systems of the form
	\begin{equation}\label{eq:dynamics_intro}
		\left\{ \begin{array}{rclr}
			\frac{\mathrm{d}}{\mathrm{d}t} x(t) & = & f(x(t)) + b( x(t) ) \alpha(t) & \text{for} \; t \in (0,\infty),\\
			x(0) & = & z
		\end{array}\right.
	\end{equation}
	where the initial value $z$ lies in an open, bounded domain $\Omega \subset \mathbb R^d$ with $C^1$ boundary, $f,b\colon \mathbb R^d\to \mathbb R^d$, and $\alpha\in \mathrm{L}^\infty_{\mathrm{loc}}(0,\infty;\mathbb R)$ is a scalar control.
	The precise assumptions on $\Omega$, $f$, and $b$ will be stated later. 
	
	Our interest lies in a novel characterization of, and approximation results for, the value function $v$ associated with the optimal control problem
	\begin{align}\label{eq:value_fct_intro}
		v(z) := \inf_{\alpha \in \mathrm{L}^\infty_{\mathrm{loc}}(0,\infty;\R)} \int_0^{\infty} \left(g(x(t)) + |\alpha(t)|^2\right)\,\mathrm{d}t \quad \text{subject to } \eqref{eq:dynamics_intro}
	\end{align}
	where $g\colon \mathbb R^d \to \mathbb R_+$. It is well known~\cite{Luk69} that if $0$ is a steady state of \eqref{eq:dynamics_intro} and $f$, $b$, and $g$ are sufficiently smooth, then $v$ can be characterized locally, i.e., in a neighborhood of the origin, as a classical solution to the Hamilton--Jacobi--Bellman (HJB) equation
	\begin{align}\label{eq:hjb_intro}
		\nabla v(z)^\top f(z) + g(z) - \tfrac{1}{4}|b(z)^\top \nabla v(z)|^2 = 0, \ \ v(0)=0.
	\end{align}
	In this case, an optimal feedback law achieving \eqref{eq:value_fct_intro} is given by 
	\begin{equation}\label{eq:opt_feedback_intro}
		\alpha_{\mathrm{opt}}(t)= u_{\mathrm{opt}}( x(t) ) := -\tfrac{1}{2}b(x(t))^\top \nabla v(x(t)).	
	\end{equation}
	Unfortunately, computing this feedback requires solving the high-dimensional nonlinear PDE \eqref{eq:hjb_intro}, a bottleneck commonly called the \emph{curse of dimensionality}. Recent years have seen renewed interest in overcoming this obstacle, and several novel computational approaches have been proposed; see, e.g., \cite{DolKK21,KalK18,OstSS22}. 
	However, to the best of our knowledge, no rigorous results on the structure or approximability of $v$ are available. By contrast, in the linear-quadratic case the value function is known to be quadratic, $v(z)= z^\top \hat{P} z$, where $\hat{P}$ solves the algebraic Riccati equation. Numerous efficient numerical methods exploit this link, often through low-rank approximation techniques; see \cite{BenS13,MasOR16,Sim16a}.
	
	In this paper, our goal is to establish a connection between $v$ and a linear infinite-dimensional system associated with a Koopman lifting of the dynamics \eqref{eq:dynamics_intro} for a specific class of feedback controls $u$. For this purpose, we extend the results from \cite{BreH23}, which show that in the uncontrolled case $\alpha\equiv 0$ with $g$ of the form $g(z)=\sum_{i=1}^\infty c_i(z)^2$, the resulting value/Lyapunov function $v$ coincides with a sum-of-squares representation $
	v(z)=\sum_{i=1}^\infty p_i(z)^2
	$ 
	described by the eigenfunctions $p_i$ of a Koopman-based operator Lyapunov equation of the form
	\begin{align}\label{eq:lyap_intro}
		\langle A \phi, \psi \rangle _P + \langle \phi, A\psi \rangle_P + \langle C\phi,C\psi \rangle_{\ell_2}
		= 0 \quad \forall \phi,\psi \in \mathcal{D}(A),
	\end{align}
	where $A$ and $C$ are determined by $f$ and $c_i$, respectively.
	To build upon the results from \cite{BreH23}, the main idea of this work can be summarized as follows. Motivated by \eqref{eq:opt_feedback_intro}, we restrict ourselves to feedback controls $\alpha(t):=u(x(t))$, where $u\colon \Omega \to \mathbb R$ is generally nonlinear. Such controls transform \eqref{eq:dynamics_intro} into the closed-loop dynamics
	\begin{align}
		\label{eq:dynamics_closed_loop}
		\left\{ \begin{array}{rclr}
			\frac{\mathrm{d}}{\mathrm{d}t} x(t) & = & f_{u}(x(t)):=f(x(t)) +b(x(t))u(x(t)) & \text{for} \; t \in (0,\infty),\\
			x(0) & = & z.
		\end{array}\right.
	\end{align}
	If the feedback $u$ is stabilizing, the associated Lyapunov function is again given by the operator Lyapunov formalism discussed in \cite{BreH23}. For this strategy to work, however, we have to rigorously analyze the Koopman operator and its preadjoint on appropriate function spaces to account for the controlled dynamics $f_u$.   
	
	\subsection{Existing literature and related results}
	The use of Koopman theory is a relatively recent development in optimal control \cite{ProBK16,KorM18} that has proven to be particularly useful in the context of data-driven engineering applications \cite{KaiKB21,KorM18,Rowetal09}. A recent overview can be found in \cite{OttR21}. 
	In \cite{VilJH21}, the authors construct an optimal feedback using the \emph{Pontryagin--Koopman operator}. They motivate this approach particularly for nonlinear systems that are not linearly stabilizable, a class we do not treat here. In contrast to their approach, which analyzes the spectral properties of the Koopman operator associated with \emph{Pontryagin's} differential equation, we define a \emph{Gramian} corresponding to the feedback and investigate its spectral properties, including nuclearity.
	For affine control inputs as in \eqref{eq:dynamics_intro}, the optimal feedback law in \cite{HuaV22} is derived as the solution of a convex optimization problem by means of the Koopman and Perron--Frobenius operators. The optimization problem is then solved with a data-driven, radial-basis-based method. While the overall goal is similar, the operator-centric framework presented here yields an operator equation that differs from the minimization-based characterization in \cite{HuaV22}. In \cite{Moyetal23}, the approach of \cite{HuaV22} is extended using a sum-of-squares technique. This is worth mentioning in our context, as we also derive a sum-of-squares representation of the value function. A summary of these results can be found in \cite{Mauetal20}. More recently, in \cite{Vai25}, a connection between the Hamilton--Jacobi equation and the Koopman operator is established by considering the Hamiltonian dynamical system associated with the optimal control problem. 
	\subsection{Contribution} 
	In this paper, we introduce a bilinear form closely resembling the solution to an algebraic Riccati equation (ARE) arising in infinite-dimensional linear-quadratic control problems \cite{Benetal07}. One of our main contributions is the investigation of the spectral properties of this bilinear form and its associated operator. For linear-quadratic control problems, similar spectral properties have been derived for solutions to operator Lyapunov equations and related quantities under various assumptions \cite{Opm10,OpmRW13,Opm20}. 
	Based on a transformation of the finite-dimensional closed-loop dynamical system \eqref{eq:dynamics_closed_loop} into a linear infinite-dimensional system acting on a weighted Sobolev space of order one, our main results are:
	\begin{enumerate}[label=(\roman*)] 
		\item In Lemma~\Rref{lemma:exponentially_stability}, we show that under regularity assumptions, the Koopman semigroup corresponding to the optimal feedback is exponentially stable.
		\item Theorem~\Rref{thm:lowrank_optimalcontrol} expresses both the value function and the optimal feedback in terms of the eigenfunctions of the value bilinear form. 
		\item In Theorem~\Rref{thm:nonlinear_op_eq}, we utilize a suitable basis transformation to derive a nonlinear operator equation on $\ell_2$ which characterizes the value bilinear form. This equation collapses to the Riccati equation in the linear-quadratic case.
	\end{enumerate}
	This article is structured as follows. Section~\Rref{sec:koopman} provides several results for the Koopman semigroup on specific function spaces that are suitable for closed-loop dynamical systems. In particular, we show that the semigroup exhibits exponential stability.
	In Section~\Rref{ssec:koopman_feedback}, we consider a standard nonlinear optimal control problem and its associated value function. For the latter, we derive a sum-of-squares representation through a value bilinear form which we relate to the optimal feedback law and a nonlinear operator equation. In Section~\Rref{sec:numerics}, we illustrate our theoretical findings by a numerical example.
	\subsection{Notation}
	For a Banach space $Z$, we denote its dual by $Z^{\ast}$ and
	write $\einner{Z,Z^{\ast}}$ for its associated dual pairing. 
	For two Banach spaces $Z_{1},Z_{2}$, let 
	$\mathcal{L}(Z_{1},Z_{2})$ be the space of bounded linear operators
	$Z_{1}\!\to\! Z_{2}$; we set $\mathcal{L}(Z):=\mathcal{L}(Z,Z)$. The space of nuclear operators is denoted by $\mathcal{N}(Z_1,Z_2)$, and we set $\mathcal{N}(Z) := \mathcal{N}(Z,Z)$. The space of self-adjoint nuclear operators over $Z$ is denoted by $\mathcal{N}_{\operatorname{sa}}(Z)$.
	If $A\colon \mathcal{D}(A)\subset Z\to Y$ is (possibly) unbounded and 
	$\mathcal{D}(A)$ is dense in $Z$, its adjoint is 
	$A^{\ast}\colon\mathcal{D}(A^{\ast})\subset Y^{\ast}\to Z^{\ast}$.
	For a bounded, open domain $\Omega\subset\mathbb{R}^{d}$ with (piecewise) $\mathcal{C}^{1}$ boundary, the Sobolev space of order $k$ and integrability $p\in[1,\infty]$ is
	$W^{k,p}(\Omega)$. The Hölder-conjugate exponent is denoted by $p^{\ast}$ and
	satisfies $1/p+1/p^{\ast}=1$. The space of smooth functions with compact support is denoted by $\CcInfty{\Omega}$. By $\clSpan \{ p_1, p_2, \dots \}$, we denote the closed span.

	To distinguish matrices from abstract operators we use  $\hat{A}$.  
	Eigenvalues of a matrix $\hat{A}$ are written $\lambda_i(\hat{A})$, and the
	Jacobian matrix of $f \colon \Omega \to \R^d$ at $z \in \Omega$ is denoted by $\mathrm{D}f(z)$. The set of all sequences over $\R$ is denoted by $\R^\N$.

	\section{Koopman and Perron--Frobenius semigroups}
	\label{sec:koopman}
	In this section, we introduce the \emph{Koopman semigroup}, which allows us to reformulate the nonlinear finite-dimensional system \eqref{eq:dynamics_closed_loop}, describing pointwise dynamics, as an infinite-dimensional linear problem and thereby to draw in part on well-known linear control theory; see, e.g., \cite{TucW09}. This so-called \emph{lifting} is well known in the literature and can be found in many textbooks, e.g., \cite{LasM94}.
	
	The Koopman operator acts on observables by evaluating them along the flow $\Phi^t(z)$ of a system of the form 
	\begin{align*}
		\label{eq:dynamics}
		\left\{ \begin{array}{rclr}
			\frac{\mathrm{d}}{\mathrm{d}t} x(t) & = & f(x(t))  & \text{for} \; t \in (0,\infty),\\
			x(0) & = & z.
		\end{array}\right.
	\end{align*}
	Formally, it is defined as an operator taking observables from a function space $Y^\ast \subseteq \{ \varphi \colon \Omega \to \mathbb{R} \}$ (for the exact definition, see Definition~\Rref{defi:weighted_sobolev}) and composing them with the closed-loop flow  
	\[  
	S^\ast(t) \colon Y^\ast \to Y^\ast, \quad \psi \mapsto \psi \circ \Phi^t  
	\]  
	for a given $t \ge 0$. In what follows, we show that, for a suitable choice of the function space $Y$ and under further assumptions on the dynamical system, the family of Koopman operators forms a weak-* continuous semigroup (see Theorem~\Rref{thm:cont_semigroup}). Let us emphasize that, for the typical choices of $Y$ and $Y^\ast$, the operators $S^\ast(t)$ are not strongly continuous \cite[Remark 5.14]{EngN06}. 
	
	Some of the following results are generalizations of our previous work \cite{BreH23} to function spaces of higher regularity (Lipschitz continuous functions). In contrast to \cite{BreH23}, we weaken the assumptions on the dynamics $f$ to incorporate the control function in the subsequent sections. Let us also note that, at a formal level, we have already written the Koopman operator as the adjoint of some semigroup $S(t)$ and identified the space of observables as a dual space. This will be made rigorous in Lemma~\Rref{lemma:preadjoint_semigroup}. 

	\subsection{Weighted Lipschitz spaces and their predual}
	\label{ssec:weighted_sobolev}
	As in \cite{BreH23}, in the following we will work with weighted Lebesgue spaces. We briefly recall these spaces before we extend some results of \cite{BreH23} to specifically weighted Sobolev spaces. This extension is needed to show higher regularity of the eigenfunctions of the value bilinear form introduced further below.

	\begin{dfntn}
		\label{defi:weighted_lp}
		Given a weight function $w\colon \Omega\to \mathbb R_+$ with $w^{-1}\in \fWp{\infty}{\Omega}$ and $p \in \{ 1,2,\infty\}$, we define
		\[
		\| \phi \|_{\Lpw{p}{w}{\Omega}} := \left\{ \begin{array}{lcl}
			\left( \int_\Omega  |\phi(x)|^p \; w(x)^p  \dint{x} \right)^{\sfrac{1}{p}}  & \qquad & \text{for} \; p \in \{ 1,2 \} \vspace{2mm}\\
			\underset{ x \in \Omega} {\esssup }\; \{ | \phi(x) | w(x) \}  & \qquad & \text{for} \; p = \infty
		\end{array}\right. 
		\]
		and the corresponding space
		\[
		\fLpw{p}{w}{\Omega} := \left\{ \phi \colon \Omega \to \R \;\mid\; \| \phi \|_{\Lpw{p}{w}{\Omega}} < \infty \right\}.
		\]
		If the domain is clear from the context, we simply write $\Lpw{p}{w}{\Omega}$ instead of $\fLpw{p}{w}{\Omega}$.
	\end{dfntn}
	The corresponding dual spaces can be characterized by \cite[Lemma 2.2]{BreH23} as follows.
	
	\begin{lmm}
		\label{lemma:dual_lpw}
		Let $p \in \{1,2\}$, where $p^\ast = \infty$ if $p=1$ and $p^\ast = 2$ if $p=2$. Then $\fLpw{p^\ast}{w}{\Omega} \simeq \left(\fLpw{p}{w}{\Omega}\right)^\ast$, with the dual pairing
		\[
		\inner[\Lpw{p}{w}{\Omega},\Lpw{p^\ast}{w}{\Omega}]{\phi}{\varphi} =
			\int_{\Omega} \phi(x) \varphi(x) w^2(x) \dint x,
		\]
		where $\phi \in \fLpw{p}{w}{\Omega}$ and $\varphi \in \fLpw{p^\ast}{w}{\Omega}$.
	\end{lmm}
	Next, we define a weighted Sobolev space of regularity one.
	\begin{dfntn}
		\label{defi:weighted_sobolev}
		Given a weight function $w\colon \Omega\to \mathbb R_+$ with $w^{-1}\in \fWp{\infty}{\Omega}$, we define the norm
		{
			\[
			\| \phi \|_{\Wpw{\infty}{w}{\Omega}} := \max  \left( \| \phi \|_{\Lpw{\infty}{w}{\Omega}}, \max_{|\alpha| = 1 } \| \mathrm{D}^\alpha \phi \|_{\Lp{\infty}{\Omega}} \right)
			\]
		}
		with the corresponding space
		\[
		\fWpw{\infty}{w}{\Omega} := \left\{ \phi \colon \Omega \to \R \; \vert \; \| \phi \|_{\Wpw{\infty}{w}{\Omega}}  < \infty  \right\}.
		\]
		For $p \in \{1,2\}$ we define the norm and the corresponding space as
		\[
		\| \phi \|_{\Wpw{p}{w}{\Omega}} :=
		\left( \| \phi \|^p_{\Lpw{p}{w}{\Omega}} + \underset{|\alpha| = 1}{\sum} \| \mathrm{D}^\alpha \phi \|^p_{\Lp{p}{\Omega}} \right)^{\sfrac{1}{p}}
		\quad \text{and} \quad
		\fWpw{p}{w}{\Omega} :=  \overline{ \fWpw{\infty}{w}{\Omega} }^{\Wpw{p}{w}{\Omega} }.
		\]
		If the domain $\Omega$ is clear from the context, we simply write $\Wpw{p}{w}{\Omega}$ instead of $\fWpw{p}{w}{\Omega}$. 
	\end{dfntn}
	\begin{rmrk}
		Note that we do not use the weighted norm for the partial derivatives, as this will not be required to show stability and continuity of the corresponding composition semigroup later on.
	\end{rmrk}
	\begin{rmrk}
		\label{rem:well_defined_weighted_sobolev}
		The space is well-defined, since for $\phi \in \fWpw{\infty}{w}{\Omega}$, we have
		\[
		\| \phi \|^p_{\Lpw{p}{w}{\Omega}} = \int_\Omega ( |\phi(x) | \, w(x) )^p \dint x \le |\Omega| \| \phi \|^p_{\Lpw{\infty}{w}{\Omega}}.
		\]
		For the partial derivatives, the boundedness of $\Omega$ yields the continuous embedding $\fLp{\infty}{\Omega}\hookrightarrow\fLp{p}{\Omega}$. Consequently, $\| \cdot\|_{\Wpw{p}{w}{\Omega}} \le C \| \cdot \|_{\Wpw{\infty}{w}{\Omega}}$, so the closure is well-defined.
	\end{rmrk}

	\begin{lmm}
		\label{lemma:density_wpwinfty}
		The space $\fWpw{\infty}{w}{\Omega}$ is a dense subset of $\fLpw{p}{w}{\Omega}$ for $p \in \{1,2\}$.
	\end{lmm}
	\begin{proof}
		Let $\phi \in \fWpw{\infty}{w}{\Omega}$. Then
		\[
		\| \phi \|^p_{\Lpw{p}{w}{\Omega}} = \int_\Omega | \phi(x) |^p w^p(x) \dint x \le \| \phi \|^p_{\Lpw{\infty}{w}{\Omega}} | \Omega| < \infty.
		\]
		Therefore, $\phi \in \fLpw{p}{w}{\Omega}$, and we conclude $\fWpw{\infty}{w}{\Omega} \subseteq \fLpw{p}{w}{\Omega}$. Next, we define the extension 
		\[
		E \colon \fLp{p}{\Omega} \to \fLp{p}{\R^{d}}, \quad \text{with} \quad E \psi := \left\{ \begin{array}{lcl}
			\psi(x) & \qquad & \text{for $x \in \Omega$}\\
			0 & \qquad & \text{else}
		\end{array}\right. .
		\]
		To show density, let $\phi \in \fLpw{p}{w}{\Omega}$ and
		\[
		\phi_k := \left( \tfrac{1}{w}  \eta_{1/k} \ast E ( \phi w ) \right) \big|_{\Omega} \in \fWpw{\infty}{w}{\Omega}
		\]
		where $\eta_\varepsilon$ denotes the standard mollifier from \cite[Chapter 4.4]{Eva98}. Here we have used $w^{-1} \in \fWp{\infty}{\Omega}$ and the smoothness of the mollified function. By \cite[Appendix, Theorem 7]{Eva98},
		\begin{align*}
			0 = \lim_{k \rightarrow \infty} \| \phi_k w - \phi w \|_{\Lp{p}{\Omega}} = \lim_{k \rightarrow \infty} \| \phi_k - \phi \|_{\Lpw{p}{w}{\Omega}}
		\end{align*}
		which was to be shown.
	\end{proof}
	Similar to $\fLp{\infty}{\Omega}$, it can be shown that $\fWpw{\infty}{w}{\Omega}$ is neither reflexive nor separable. Moreover, its dual space does not admit a simple or explicit characterization. For this reason, its predual plays a central role in the subsequent analysis. We therefore begin by introducing an appropriate space and later demonstrate that this space can be used to identify the predual of $\fWpw{\infty}{w}{\Omega}$.
	\begin{dfntn}
		\label{defi:neg_weighted_sobolev}
		Let $p \in \{1,2\}$, where $p^\ast = \infty$ if $p=1$ and $p^\ast = 2$ if $p=2$. We define the space
		\[
		\fNWpw{p}{w}{\Omega} := \left\{ \; [(f_0,f_1, \dots, f_d)]_{\sim}  \; \mid \; f_0 \in \Lpw{p}{w}{\Omega} \; \text{and} \; f_1, \dots, f_d \in \Lp{p}{\Omega} \right\} 
		\]
		with equivalence relation
		\begin{align*}
			& (f_0, f_1, \dots, f_d ) \sim (\tilde{f}_0, \tilde{f}_1, \dots, \tilde{f}_d) \\
			\Leftrightarrow \quad  & \inner[\Lpw{p}{w}{\Omega},\Lpw{p^\ast}{w}{\Omega}]{f_0 - \tilde{f}_0}{g} + \sum_{i=1}^d \inner[\Lp{p}{\Omega},\Lp{p^\ast}{\Omega}]{f_i - \tilde{f}_i }{ \partial_{i} g} = 0 \qquad \forall g \in \fWpw{p^\ast}{w}{\Omega}.	
		\end{align*}
		Furthermore, we endow the space with the norm
		\[
		\| [f]_{\sim} \|^p_{\fNWpw{p}{w}{\Omega}} := \underset{ \tilde{f} \sim f }{\inf} \; \left\{ \| \tilde{f}_0 \|^p_{\Lpw{p}{w}{\Omega}} + \sum_{i=1}^d \| \tilde{f}_i \|^p_{\Lp{p}{\Omega}} \right\} .
		\]
	\end{dfntn}
	As is commonly done for classical Lebesgue spaces, we identify an equivalence class with an arbitrarily chosen representative whenever convenient.
	\begin{rmrk}
		\label{rem:invariance_equiv_rel}
		Note that the equivalence relation remains the same across spaces. 
		In particular, for any representative $f = (f_0, \dots, f_{d})$ with $f_0 \in  \Lpw{2}{w}{\Omega}$ and $f_1, \dots, f_d \in \Lp{2}{\Omega}$ 
		the equivalence class $[f]_{\sim} \in \fNWpw{2}{w}{\Omega}$ can be naturally extended to $[f]_{\sim} \in \fNWpw{1}{w}{\Omega}$. 
		This follows directly from the density of $\fWpw{\infty}{w}{\Omega}$ in both spaces, as well as the fact that the dual pairing remains unchanged. Thus a common representative determines the same equivalence class in either space.
	\end{rmrk}
	\begin{rmrk}
		Note that the definition differs from the usual definition of $H^{-1}(\Omega)$ as $H_0^{1}(\Omega)^\ast$ found in the literature \cite[Chapter 5.9]{Eva98}. Our aim is to obtain a characterization of the predual of the Banach space $\fWpw{\infty}{w}{\Omega}$ analogous to the identification $H^1(\Omega)^\ast \simeq H^1(\Omega)$ provided by the Riesz representation theorem. We subsequently pass to the Hilbert space $\fWpw{2}{w}{\Omega}$ by an embedding.
	\end{rmrk}
	\begin{lmm}
		\label{lemma:dual_Wpw}
		Let $p \in \{1,2\}$, where $p^\ast = \infty$ if $p=1$ and $p^\ast = 2$ if $p=2$. Then
		\[
		(\fNWpw{p}{w}{\Omega})^\ast \simeq \fWpw{p^\ast}{w}{\Omega} 
		\]
		where the dual pairing between $\phi \in \fNWpw{p}{w}{\Omega}$ and $\psi \in \fWpw{p^\ast}{w}{\Omega}$ is defined as 
		\[
		\inner{\phi}{\psi} := \inner[\Lpw{p}{w}{\Omega}, \Lpw{p^\ast}{w}{\Omega}]{ \phi_0 }{\psi} + \sum_{i=1}^d  \inner[\Lp{p}{\Omega},\Lp{p^\ast}{\Omega}]{ \phi_i }{ \partial_{i} \psi }. 
		\]
	\end{lmm}
	\begin{proof}
		For the proof, we follow the notation from \cite[Chapter 4]{Rud91} and define
		\[
		Z := \fLpw{p}{w}{\Omega} \times \fLp{p}{\Omega}^d \quad \text{and} \quad Z^\ast := \fLpw{p^\ast}{w}{\Omega} \times \fLp{p^\ast}{\Omega}^d.
		\]
		Indeed, with \cite[Theorem 6.12]{Alt12} and Lemma~\Rref{lemma:dual_lpw}, we conclude that $Z^\ast$ is the dual space of $Z$ with respect to the pairing
		\[
		\inner[Z,Z^\ast]{f}{g} := \inner[\Lpw{p}{w}{\Omega},\Lpw{p^\ast}{w}{\Omega}]{f_0}{g_0} + \sum_{i=1}^d \inner[\Lp{p}{\Omega},\Lp{p^\ast}{\Omega}]{f_i}{g_i}.
		\]
		Next, we consider the embedding
		\[
		\Psi \colon  \fWpw{p^\ast}{w}{\Omega}  \to Z^\ast, \qquad 
		 \phi  \mapsto  ( \phi, \partial_{1} \phi, \dots, \partial_{d} \phi ).
		\]
		It is easy to see that $\Psi$ is isometric. We define 
		\[
		N := \Range \Psi \quad \text{and} \quad M := \prescript{\perp}{}{N} = \{ \phi \in Z \; \big| \; \inner[Z,Z^\ast]{\phi}{\psi} = 0 \quad \text{for all} \; \psi \in N \}.
		\]
		As mentioned in \cite[Section 4.6]{Rud91}, $M$ is closed with respect to the norm topology. Let us show that $(\prescript{\perp}{}{N})^\perp = N$. First, note that $N \subseteq (\prescript{\perp}{}{N})^\perp$ by \cite[Theorem 4.7]{Rud91}.
		
		To show that $ (\prescript{\perp}{}{N})^\perp \subseteq N$, we first consider an arbitrary
		$\phi \in \CcInfty{\Omega}$ and $1 \le i \le d$. Integration by parts gives
		\[
		\int_\Omega \left( \tfrac{\partial_i \phi}{w^2} \right) g \; w^2 \dint x + \int_\Omega \partial_i g \, \phi \dint x = 0 \qquad \text{for all}\; g \in \fWpw{\infty}{w}{\Omega}.
		\]
		Since, by assumption, $ \sfrac{1}{w} \in W^{1,\infty}(\Omega)$ and $\Omega$ is bounded, we conclude that
		\[
		\| \tfrac{\partial_i \phi}{w^2}\|_{\Lpw{p}{w}{\Omega}} = \left( \int_\Omega \tfrac{|\partial_i \phi|^p}{w^{2 p}}  w^p \dint x \right)^{\sfrac{1}{p}}   \le |\Omega|^{\sfrac{1}{p}} \; \| \tfrac{1}{w} \|_{\Lp{\infty}{\Omega}} \|  \partial_i \phi \|_{\Lp{\infty}{\Omega}} < \infty,
		\]
		which implies that
		\begin{equation}
			\label{eq:test_fct}
			\left( \left( \sfrac{\partial_i \phi}{w^2} \right), 0,\dots, 0, \phi, 0, \dots, 0 \right) \in \prescript{\perp}{}{N}.
		\end{equation}
		We are ready to show that $ (\prescript{\perp}{}{N})^\perp \subseteq N$. To this end, let $\psi \in (\prescript{\perp}{}{N})^\perp$. With \eqref{eq:test_fct}, we obtain 
		\[
		0 = \int_\Omega \psi_0 \, \partial_i \phi \dint x +  \int_\Omega \psi_i \, \phi \dint x  \qquad \text{for all} \; \phi \in \CcInfty{\Omega}.
		\]
		However, this implies that $\psi_i$ is the weak derivative of $\psi_0$ by definition and therefore $\psi \in N$.
		All assumptions of \cite[Theorem 4.9 (b)]{Rud91} are fulfilled and we conclude
		\[
		\fWpw{p^\ast}{w}{\Omega} \simeq N = M^\perp \simeq (Z/M)^\ast = \fNWpw{p}{w}{\Omega}^\ast.
		\]
	\end{proof}
	\begin{rmrk}
		\label{rem:arens_eels}
		For the weight $w(x) = \sfrac{1}{\|x\|}$, Rademacher's theorem \cite[Theorem 4.5]{EvaG92} allows us to identify $\fWpw{\infty}{w}{\Omega}$ as the space of Lipschitz continuous functions that vanish at $0$. This space is often denoted by $\operatorname{Lip}_0(\Omega)$. It is well known, e.g., from \cite{Wea99}, that $\operatorname{Lip}_0(\Omega)$ is a dual space and its predual can be characterized \cite[Chapter 3]{Wea99} by so-called molecules. In more detail, let us define the \emph{Arens--Eells space}
		\begin{align*}
			& \text{\AE}(\Omega) := \operatorname{cl} \; {\left\{ \sum_{i=1}^n a_i \left( \delta_{p_i} - \delta_{q_i} \right)\quad \text{where} \quad p_i,q_i \in \Omega \; \text{and} \; a_i \in \R \right\}} 	\\
			\text{with respect to} \qquad & \| m \|_{\text{\AE}(\Omega)} := \inf \left\{ \sum_{i=1}^n | a_i | \| p_i - q_i \|_{\R^d} \quad \text{such that} \quad m= \sum_{i=1}^n a_i \left( \delta_{p_i} - \delta_{q_i} \right) \right\}.
		\end{align*}
		Then $\operatorname{Lip}_0(\Omega) \simeq \text{\AE}(\Omega)^\ast$ \cite[Theorem 3.3]{Wea99}, where the identification is made via
		\[
		\inner[\text{\AE}(\Omega),\operatorname{Lip}_0(\Omega)]{\sum_{i=1}^n a_i \left( \delta_{p_i} - \delta_{q_i}\right)}{f} := \sum_{i=1}^n a_i \left( f(p_i) - f(q_i) \right).
		\]
		Remarkably, this characterization generalizes to infinite-dimensional state spaces, i.e., it is possible to consider $\Omega \subseteq X$ where $X$ is a generic metric space. However, we do not use this identification here for several reasons, including the need for an embedded Hilbert space, as shown in Lemma~\Rref{lemma:gelfand_like}.
	\end{rmrk}
	As in the usual Sobolev case with $p=2$, the weighted space $\fWpw{2}{w}{\Omega}$ has a Hilbert space structure, which allows us to identify its predual with the space itself via the Riesz representation theorem.
	\begin{lmm}
		\label{lemma:iso_hilbert}
		We have
		\[
		\fWpw{2}{w}{\Omega} \simeq \fNWpw{2}{w}{\Omega},
		\]
		with the isometric isomorphism
		\[
		\Psi \colon  \fWpw{2}{w}{\Omega}  \to \fNWpw{2}{w}{\Omega}, \qquad
		\phi  \mapsto  [ ( \phi, \partial_{1} \phi, \dots, \partial_{d} \phi)]_{\sim}.
		\]
	\end{lmm}
	\begin{proof}
		On $\fWpw{2}{w}{\Omega}$ we define the following norm-inducing scalar product
		\[
		\inner[\Wpw{2}{w}{\Omega}]{\phi}{\psi} := \int_\Omega \phi(x) \psi(x) w^2(x) \dint x + \sum_{i = 1}^d \int_\Omega \partial_{i} \phi(x) \partial_{i} \psi(x) \dint x.
		\]
		By the Riesz representation theorem \cite[Theorem 6.1]{Alt12} and Lemma~\Rref{lemma:dual_Wpw}, for each
		\[
		f \in \fNWpw{2}{w}{\Omega} \simeq \left( \fWpw{2}{w}{\Omega} \right)^\ast \simeq \fWpw{2}{w}{\Omega}
		\]
		there exists a unique $\psi_f \in \fWpw{2}{w}{\Omega}$ such that
		\[
		\inner[\Wpw{2}{w}{\Omega}]{\phi}{\psi_f}  = \inner[\Wpw{2}{w}{\Omega},\NWpw{2}{w}{\Omega}]{\phi}{f} 	\qquad \text{for all}\; \phi \in \fWpw{2}{w}{\Omega},
		\]
		which further satisfies $\| \psi_f \|_{\Wpw{2}{w}{\Omega}} = \| f \|_{\NWpw{2}{w}{\Omega}}$. Therefore $\Psi$ is bijective and
		\[
		( \psi_f, \partial_{1} \psi_f, \dots, \partial_{d} \psi_f) \in f \quad \text{with} \quad \| f \|_{\NWpw{2}{w}{\Omega}} = \| \psi_f \|_{\Wpw{2}{w}{\Omega}}.
		\]
	\end{proof}
	Hence, these spaces have a Gelfand-like structure in which the roles of the space and its dual are reversed. 
	\begin{lmm}
		\label{lemma:gelfand_like}
		We have
		\[
		\fWpw{\infty}{w}{\Omega} \subset \fWpw{2}{w}{\Omega} \simeq	\fNWpw{2}{w}{\Omega} \subset \fNWpw{1}{w}{\Omega},
		\]
		where the embeddings are dense.
	\end{lmm}
	\begin{proof}
		The density of the embedding $\fWpw{\infty}{w}{\Omega} \subset \fWpw{2}{w}{\Omega}$ follows directly from Definition~\Rref{defi:weighted_sobolev} and Remark~\Rref{rem:well_defined_weighted_sobolev}. The isomorphism $\fWpw{2}{w}{\Omega} \simeq	\fNWpw{2}{w}{\Omega}$ was already shown in Lemma~\Rref{lemma:iso_hilbert}. It remains to show that $\fNWpw{2}{w}{\Omega} \subset \fNWpw{1}{w}{\Omega}$ is dense. 
		
		To this end, let $ [ ( \phi_0, \dots,\phi_d )]_{\sim} \in \fNWpw{1}{w}{\Omega} $ have the arbitrarily chosen representative $(\phi_0, \dots,\phi_d)$. By \cite[Chapter 5.3.3, Theorem 3]{Eva98}, $\fLp{2}{\Omega}$ is dense in $\fLp{1}{\Omega}$, and by \cite[Lemma 2.3]{BreH23}
		$\fLpw{2}{w}{\Omega}$ is dense in $\fLpw{1}{w}{\Omega}$. Thus, there exist
		\begin{align*}
			& \phi_0^{(n)} \in  \fLpw{2}{w}{\Omega} \subset \fLpw{1}{w}{\Omega} \quad \text{such that}\quad \lim_{n \rightarrow  \infty} \| \phi_0^{(n)} - \phi_0 \|_{\Lpw{1}{w}{\Omega}} = 0 \\
			\text{and} \qquad & \phi_i^{(n)} \in  \fLp{2}{\Omega} \subset \fLp{1}{\Omega} \quad \text{such that}\quad \lim_{n \rightarrow  \infty} \| \phi_i^{(n)} - \phi_i \|_{\Lp{1}{\Omega}} = 0,
		\end{align*}
		for all $1 \le i \le d$. Take $\phi_n = [ ( \phi_0^{(n)}, \dots,\phi_d^{(n)} )]_{\sim} \in \fNWpw{2}{w}{\Omega}$. As already stated in Remark~\Rref{rem:invariance_equiv_rel}, the equivalence relation does not change between spaces, and therefore the choice of the representative is not relevant. 
		Furthermore, 
		\begin{align*}
			& \lim_{n \rightarrow \infty} \| \phi - \phi_n \|_{\NWpw{1}{w}{\Omega}} 
			= \lim_{n \rightarrow \infty} \underset{ f \in [\phi - \phi_n]_{\sim}}{\inf} \;  \left\{ \| f_0 \|_{\Lpw{1}{w}{\Omega}} + \sum_{i=1}^d \| f_i \|_{\Lp{1}{\Omega}} \right\} \\
			\le & \lim_{n \rightarrow \infty} \left( \| \phi_0 - \phi_0^{(n)} \|_{\Lpw{1}{w}{\Omega}} + \sum_{i=1}^d \| \phi_i - \phi_i^{(n)} \|_{\Lp{1}{\Omega}} \right) = 0.
		\end{align*}
	\end{proof}
	Analogously to Lemma~\Rref{lemma:density_wpwinfty}, we can show that there exists a more regular dense subset of $\fNWpw{1}{w}{\Omega}$, which enables us to first establish certain properties on this subset and then extend them to the whole space by a density argument.
	\begin{lmm}
		\label{lemma:dense_zeros}
		The set 
		\[	
		D := \left\{
		[(\phi_0, 0, \dots, 0) ]_{\sim}
		\;\mid\;
		\phi_0 \in \fWpw{\infty}{w}{\Omega}
		\right\}
		\]
		is dense in $\fNWpw{1}{w}{\Omega}$ and $\fNWpw{2}{w}{\Omega}$.
	\end{lmm}
	\begin{proof}
		Let $p \in \{1,2\}$, where $p^\ast = \infty$ if $p=1$ and $p^\ast = 2$ if $p=2$. Let $[f]_{\sim} \in \fNWpw{p}{w}{\Omega}$. By definition, there exists a representative
		\[
		(f_0, \dots, f_d ) \in [f]_{\sim}  \quad \text{with}\; f_0 \in \fLpw{p}{w}{\Omega} \; \text{and}\; f_i \in \fLp{p}{\Omega} \; \text{for}\; 1 \le i \le d.
		\]
		As already shown in Lemma~\Rref{lemma:density_wpwinfty}, under the stated assumptions $\fWpw{\infty}{w}{\Omega}$ is dense in $\fLpw{p}{w}{\Omega}$; furthermore, $\CcInfty{\Omega}$ is dense in $\fLp{p}{\Omega}$~\cite[Appendix C, Theorem 6]{Eva98}. 
		Therefore, there exist $\bar{f}_0 \in \fWpw{\infty}{w}{\Omega}$ and $\bar{f}_i \in \CcInfty{\Omega}$ such that
		\[
		\| f_0 - \bar{f}_0 \|_{\Lpw{p}{w}{\Omega}} \le \varepsilon \quad \text{and} \quad \| f_i - \bar{f}_i \|_{\Lp{p}{\Omega}} \le \varepsilon \quad \text{for }1\le i\le d.
		\]
		Thus we obtain
		\begin{align}
			\label{eq:dense_zeros:1}
		\| [f]_{\sim} - [\bar{f}]_{\sim} \|_{\NWpw{p}{w}{\Omega}}^p = & \inf_{(g_0,\dots,g_d) \in [f - \bar{f}]_{\sim}} \; \left\{ \| g_0 \|_{\Lpw{p}{w}{\Omega}}^p + \sum_{i=1}^d \| g_i \|_{\Lp{p}{\Omega}}^p \right\} \\
		\le &   \| f_0 - \bar{f}_0 \|_{\Lpw{p}{w}{\Omega}}^p + \sum_{i=1}^d \| f_i - \bar{f}_i \|_{\Lp{p}{\Omega}}^p  \le (d+1)\varepsilon^p.
		\end{align}
		However, for all $\phi \in \fWpw{\infty}{w}{\Omega}$, we have
		\begin{align*}
			\inner[\NWpw{p}{w}{\Omega},\Wpw{p^\ast}{w}{\Omega}]{\bar{f}}{\phi} = & \inner[\Lpw{p}{w}{\Omega},\Lpw{p^\ast}{w}{\Omega}]{\bar{f}_0}{\phi} + \sum_{i=1}^d \inner[\Lp{p}{\Omega},\Lp{p^\ast}{\Omega}]{\bar{f}_i}{ \partial_i \phi}\\
			= & \inner[\Lpw{p}{w}{\Omega},\Lpw{p^\ast}{w}{\Omega}]{\bar{f}_0}{\phi} - \sum_{i=1}^d \inner[\Lp{p}{\Omega},\Lp{p^\ast}{\Omega}]{\partial_i \bar{f}_i}{\phi}
			= : \inner[\Lpw{p}{w}{\Omega},\Lpw{p^\ast}{w}{\Omega}]{g}{\phi},
		\end{align*}
		where 
		$
		g := 
		\bar{f}_0- \sum_{i=1}^d \frac{1}{w^2}(\partial_i \bar{f}_i)  .
		$
		Since $\sfrac{1}{w} \in W^{1,\infty}(\Omega)$ by assumption, we conclude 
		\[
		\|\sfrac{\partial_i \bar{f}_i}{w^2} \|_{\Lpw{\infty}{w}{\Omega}} \le \| \sfrac{1}{w} \|_{\Lp{\infty}{\Omega}} \| \partial_i \bar{f}_i \|_{\Lp{\infty}{\Omega}} < \infty.
		\]
		This means $g \in \fLpw{\infty}{w}{\Omega}$. Since $\Omega$ is bounded, $w^{-2}$ is also Lipschitz continuous and we conclude that $g \in \fWpw{\infty}{w}{\Omega}$. 
		In particular, $g \in \fLpw{p}{w}{\Omega}$, and we further conclude that $( g, 0, \dots, 0) \in [\bar{f}]_{\sim}$. With \eqref{eq:dense_zeros:1}, it follows that
		\[
		\| [f]_{\sim} - [(g,0,\dots,0)]_{\sim} \|_{\NWpw{p}{w}{\Omega}} = \| [f]_{\sim} - [\bar{f}]_{\sim} \|_{\NWpw{p}{w}{\Omega}} \le (d+1)^{1/p}\varepsilon,
		\]
		which shows that any $[f]_{\sim}$ can be approximated arbitrarily well by elements of $D$.
	\end{proof}
	\begin{rmrk}
		Let us emphasize that we cannot express the norm of the equivalence class by the $\fLpw{p}{w}{\Omega}$ norm of the representative, i.e.,
		\[
		\|[(\phi_0,0,\dots,0)]_{\sim} \|_{\NWpw{1}{w}{\Omega}} \neq \| \phi_0 \|_{\Lpw{1}{w}{\Omega}}.
		\]
		However, the result allows us to see the space $\fNWpw{1}{w}{\Omega}$ as the completion of the space $\fLpw{1}{w}{\Omega}$ with a weaker norm.
	\end{rmrk}
	\subsection{The Koopman semigroup on $\fWpw{\infty}{w}{\Omega}$}
	\label{ssec:koopman_uncontrolled}
	
	Throughout the remainder of this manuscript, we fix $w \colon \Omega \to \R_+$ with $w^{-1} \in \fWp{\infty}{\Omega}$ and adopt the following notation:
	\[
	Y := \fNWpw{1}{w}{\Omega}, \quad 
	Y^\ast := \fWpw{\infty}{w}{\Omega}, \quad 
	X := \fLpw{1}{w}{\Omega}, \quad 
	X^\ast := \fLpw{\infty}{w}{\Omega}.
	\]
	The corresponding Hilbert spaces are defined as 
	\[
	H_X := \fLpw{2}{w}{\Omega}, \quad 
	H_Y := \fWpw{2}{w}{\Omega}.
	\]
	In \cite{BreH23}, the Koopman semigroup and the Perron--Frobenius semigroup were analyzed on the spaces $X^\ast$ and $X$ under the assumption $f \in \mathcal{C}^{1}(\Omega)$. However, when considering a suitable Koopman lifting for the closed-loop dynamics, it is more appropriate to work in the spaces $Y^\ast$ and $Y$ and to relax the regularity assumptions on $f$. For $f \in \fWpw{\infty}{w}{\Omega}$,  
	this allows us to consider the feedback law as an element of the space on which the semigroup operates. 
	
	We consider dynamical systems with inward-pointing boundary dynamics to ensure the well-posedness of the Koopman semigroup introduced later. We therefore make the following definitions.
	\begin{dfntn}
		\label{defi:bd_admissible_dynamic}
		We call $f \colon \Omega \to \R^d$ a \emph{boundary-admissible dynamic} if
		\begin{enumerate}[label=(\roman*)]
			\item $ f_i \in \fWpw{\infty}{w}{\Omega} $ for $1 \le i \le d$.
			\label{defi:bd_admissible_dynamic_regularity}
			\item $ \nu(x)^\top f(x) \le 0 $ for all $ x \in \partial \Omega $, where $\nu$ denotes the outward unit normal.
			\label{defi:bd_admissible_dynamic_tangent}
		\end{enumerate}
	\end{dfntn}
	\begin{dfntn}
		\label{defi:uncontrolled_flow}
		Let $f$ be a boundary-admissible dynamic. The \emph{flow of $f$} is the mapping $\Phi^t \colon \Omega \to \Omega$ such that for every $z \in \Omega$ 
		\begin{equation}
			\left\{ 
			\begin{array}{rclcl}
				\tfrac{\mathrm{d}}{\mathrm{d} t} \Phi^t(z) & = & f( \Phi^t(z)) & \qquad & \text{for } 0 < t < \infty\\
				\Phi^0(z) & = & z & \qquad & 
			\end{array}
			\right.
		\end{equation}
	\end{dfntn}
	We proceed to justify the well-posedness of the flow.
	\begin{prpstn}
		\label{prop:bounds_flow}
		Let $f$ be a boundary-admissible dynamic. Then its flow $\Phi^t \colon \Omega \to \Phi^t(\Omega)$ is well-defined and Lipschitz continuous with constant
		\[
		L := \exp( \| f \|_{\Wpw{\infty}{w}{\Omega}} t ).
		\]
		Furthermore, the inverse $ \Phi^{-t} \colon \Phi^t(\Omega) \to \Omega$ exists and is Lipschitz continuous with the same constant, and
		\[
		0 < L^{-d} \le \bigl|\det{\mathrm{D} \Phi^t(z)}\bigr| \le L^{d} < \infty \qquad \text{for almost every}\; z \in \Omega.
		\]
	\end{prpstn}
	\begin{proof}
		See Appendix~\Rref{sec:technical_proofs}.
	\end{proof}
	\begin{prpstn}
		\label{prop:explicit_Dphi}
		Let $f$ be a boundary-admissible dynamic, let $\Phi^t$ be its flow, and let $T >0$. Then, for almost every $z\in \Omega$, the map $t \mapsto \mathrm{D} \Phi^t(z)$ is the unique solution in the Carath\'{e}odory sense of the linear matrix-valued ODE 
		\[
		\left\{
		\begin{array}{rclrl}
		\tfrac{\mathrm{d}}{\mathrm{d}t} \mathrm{D} \Phi^t(z) &=& \mathrm{D}f( \Phi^t(z)) \mathrm{D}\Phi^t(z)& \qquad \text{for almost every}\; & t \in (0,T)\\
			\mathrm{D} \Phi^0(z) & =& I.
		\end{array}
		\right.
		\]
	\end{prpstn}
	\begin{proof}
		See Appendix~\Rref{sec:technical_proofs}.
	\end{proof}
	\begin{dfntn}
		\label{defi:composition_semigroup}
		Let $Z \in \{ X, H_X, Y, H_Y \}$ and $f$ be a boundary-admissible dynamic with flow $\Phi^t$. For $t \in [0,\infty)$, we define the Koopman semigroup with respect to $f$ as
		\[
		S^\ast(t) \colon Z^\ast \to Z^\ast, \quad \phi \mapsto \phi \circ \Phi^t.
		\]
	\end{dfntn}
	To establish that the Koopman semigroup introduced in Definition~\Rref{defi:composition_semigroup} has the properties of a weak-* continuous semigroup, we begin with a series of propositions. As a first step, we examine its well-posedness and boundedness.
	\begin{lmm}
		\label{lemma:dissipation_type}
		Let $f$ be a boundary-admissible dynamic with flow $\Phi^t$, and let $w$ denote the weight function from Definition~\Rref{defi:weighted_sobolev}. Then
		\[
		w(z) \le \exp(t \omega_0 ) w(\Phi^t(z)) \qquad \text{for} \quad t \in [0,\infty), \; z \in \Omega,
		\]
		where
		\[
		\omega_0 := \esssup_{z \in \Omega} \; \left\{-\frac{f(z)^\top \nabla w(z)}{w(z)}\right\} < \infty.
		\]
	\end{lmm}
	\begin{proof}
		By assumption, $f \in \fLpw{\infty}{w}{\Omega}$ and $\tfrac{1}{w} \in \fWp{\infty}{\Omega}$. A simple calculation reveals, for almost every $z \in \Omega$, that
		\begin{align*}
			- \frac{f(z)^\top \nabla w(z)}{ w(z)} =  w(z) f(z)^\top \nabla \left( \frac{1}{w(z)} \right)
			\le \| f \|_{\Lpw{\infty}{w}{\Omega}} \| \tfrac{1}{w} \|_{\Wp{\infty}{\Omega}} < \infty.
		\end{align*}
		For $t\ge s\ge0$, differentiation of $w$ along trajectories yields 
		\begin{align*}
			\frac{\mathrm{d}}{\mathrm{d}s} w(\Phi^{t-s}(z) ) = - f(\Phi^{t-s}(z))^\top \nabla{w(\Phi^{t-s}(z))} \le \omega_0 w(\Phi^{t-s}(z)).
		\end{align*} 
		The assertion now follows from Gronwall's lemma~\cite[§29, Theorem VI]{Wal98}.
	\end{proof}
	\begin{dfntn}
		\label{defi:preadjoint_semigroup}
		Let $f$ be a boundary-admissible dynamic and
		\[
		\begin{array}{rlcl}
			S(t) \colon & Y & \to & Y, \\
			& [(\rho_0, \rho_1, \dots, \rho_d)]_{\sim} & \mapsto & [(\xi_0, \xi_1, \dots, \xi_d )]_{\sim},
		\end{array}
		\]
		where for almost every $a \in \Phi^t(\Omega)$ we set
		\[
		\left\{ 
		\begin{array}{lcl}
			\xi_0(a) &=& \mu_t(a) \rho_0( \Phi^{-t}(a)) \\
			\left( \xi_1(a), \dots, \xi_d(a) \right)^\top  &=&  M_t(a) \left( \rho_1( \Phi^{-t}(a)), \dots, \rho_d( \Phi^{-t}(a)) \right)^\top			
		\end{array}
		\right.	
		\]	
		and $\xi_i\big|_{\Omega \setminus \Phi^t(\Omega)} \equiv 0$ for all $0 \le i \le d$. Furthermore, $\mu_t$ and $M_t$ are defined as \vspace{-5px}
		\begin{align*}
			\mu_t(a) &:= | \det{ \mathrm{D}\Phi^t( \Phi^{-t}(a)) } |^{-1} \; \frac{w^2( \Phi^{-t}(a))}{w^2(a)} \ge 0 \\
			\text{and}\qquad M_t(a) &:=  | \det{ \mathrm{D}\Phi^t(\Phi^{-t}(a)) } |^{-1}  \mathrm{D} \Phi^{t}( \Phi^{-t}(a))
		\end{align*}
		for $a \in \Phi^t(\Omega)$. We call $S(\cdot)$ the \emph{Perron--Frobenius semigroup} with respect to $f$ acting on $Y$. Similarly,
		\[
			S(t) \colon X \to X, \qquad \rho_0 \mapsto \xi_0
		\]
		is the \emph{Perron--Frobenius semigroup} with respect to $f$ acting on $X$.
	\end{dfntn}
	\begin{lmm}
		\label{lemma:preadjoint_semigroup}
		Let $Z \in \{X,H_X,Y,H_Y\}$ and $f$ be a boundary-admissible dynamic. The adjoint $(S(t))^\ast \colon Z^\ast \to Z^\ast$ of the Perron--Frobenius semigroup associated with $f$ is the Koopman semigroup from Definition~\Rref{defi:composition_semigroup} associated with $f$.
	\end{lmm}
	\begin{proof}
		Let $\psi \in Y$ and $\phi \in Y^\ast$. Then, for any $t \in (0,\infty)$, we have
		\begin{align*}
			\inner[Y,Y^\ast]{\psi}{S^\ast(t) \phi} = & \int_\Omega \psi_0(x) \phi( \Phi^t(x)) w^2(x) \dint x \\
			& \qquad  + \int_\Omega \left( \psi_1(x),\dots, \psi_d(x) \right)^\top \nabla \left( \phi \circ \Phi^t \right)(x) \dint x. \\
			\intertext{By the area rule \cite[Theorem 3.8]{EvaG92} and the chain rule \cite[Theorem 9.15]{Rud76} we obtain} 
			\inner[Y,Y^\ast]{\psi}{S^\ast(t) \phi} = & \int_{\Phi^t(\Omega)} \psi_0( \Phi^{-t}(x)) \phi(x) \mu_t(x) w^2(x) \dint x \\
			& \qquad + \int_{\Phi^t(\Omega)} \left( \psi_1(\Phi^{-t}(x)), \dots, \psi_d(\Phi^{-t}(x)) \right)^\top M_t(x)^\top \nabla \phi(x) \dint x.
		\end{align*}
		The case $Z = X$ follows by considering only the first term. The case $Z \in \{H_X,H_Y\}$ is analogous since the dual pairing does not change.
	\end{proof}
	\begin{prpstn}
		\label{prop:bound_semigroup}
		Let $f$ be a boundary-admissible dynamic with corresponding Koopman semigroup $S^\ast$, let $T < \infty$, and let $Z \in \{ X, H_X, Y, H_Y \}$. Then there exists some $M < \infty$ such that
		\[
		\| S^\ast(t) \|_{\mathcal{L}(Z^\ast)} \le M \qquad \text{for all}\; t \in [0,T],
		\] 
		where $M$ depends only on $\| f \|_{\Wpw{\infty}{w}{\Omega}}$, $\|\sfrac{1}{w} \|_{\Wp{\infty}{\Omega}}$, $T$, and $ \Omega$.
	\end{prpstn}
	\begin{proof}
		We begin with $Z = X = \fLpw{1}{w}{\Omega}$. With Lemma~\Rref{lemma:dissipation_type}, we derive
		\begin{align*}
			\| S^\ast(t) \phi \|_{X^\ast} = & \esssup_{z \in \Omega} | \phi(\Phi^t(z)) | w(z) \\
			\le & \esssup_{z \in \Omega} \tfrac{w(z)}{w(\Phi^t(z))} \esssup_{z \in \Omega} | \phi(\Phi^t(z)) | w(\Phi^t(z)) \le \exp( t \omega_0) \| \phi \|_{X^\ast}
		\end{align*}
		and therefore $\| S^\ast(t) \|_{\mathcal{L}(X^\ast)} \le \exp(\omega_0 t)$. Next, we consider $Z = H_X = \fLpw{2}{w}{\Omega}$. The identity from Lemma~\Rref{lemma:preadjoint_semigroup} leads to
		\[
		 \inner[H_X]{S(t) S^\ast(t) \phi}{\phi} = \int_{\Phi^t(\Omega)} | \det{ \mathrm{D}\Phi^t( \Phi^{-t}(a)) } |^{-1} \; \tfrac{w^2(\Phi^{-t}(a))}{w^2(a)} \phi^2(a) w^2(a) \dint a.
		\]
		From Lemma~\Rref{lemma:dissipation_type}, we derive $\tfrac{w^2(\Phi^{-t}(a))}{w^2(a)} \le \exp(2\omega_0 t)$ for all $a \in \Phi^t(\Omega)$. With Proposition~\Rref{prop:bounds_flow}, we then conclude
		 \begin{equation*}
		 	\| S^\ast(t) \phi \|_{H_X}^2 = \inner[H_X]{S(t) S^\ast(t) \phi}{\phi}  \le L^{d} \exp(2\omega_0 t )\| \phi \|_{H_X}^2.
		 \end{equation*}
		For $Z \in \{Y, H_Y\}$, we first note that, by Rademacher's theorem~\cite[Theorem 4.5]{EvaG92} and the Lipschitz continuity of $\Phi^t$, we can bound
		 \begin{equation}
		 	  \esssup_{z \in \Omega}\|\mathrm{D} \Phi^t(z)\|_p \le C_p < \infty \qquad \text{for}\; p \in [1,\infty) \cup \{ \infty \}.
		 	  \label{eq:prop:bound_semigroup:boundDphi}
		 \end{equation}
		For $Z = Y$ we can use this inequality and the chain rule \cite[Theorem 9.15]{Rud76} to derive  
		\begin{align*}
			\| \nabla ( S^\ast(t) \phi )(z) \|_\infty = & \| \mathrm{D}\Phi^t(z)^\top \nabla \phi( \Phi^t( z )) \|_\infty \\ \le & \| \mathrm{D}\Phi^t(z)^\top \|_\infty \| (\nabla \phi)( \Phi^t( z ))\|_\infty \le C_{\infty}  \| \phi \|_{\Wpw{\infty}{w}{\Omega}}.
		\end{align*}
		Since $\|S^\ast(t)\|_{\mathcal{L}(X^\ast)} \le \exp(\omega_0 t)$, this then leads to
		\[
		\| S^\ast(t) \phi \|_{Y^\ast} = \max \left( \| S^\ast(t) \phi\|_{\Lpw{\infty}{w}{\Omega}}, \max_{| \alpha | = 1}\| \mathrm{D}^\alpha (S^\ast(t) \phi ) \|_{\Lp{\infty}{\Omega}}  \right) \le \left( \exp(\omega_0 t) + C_{\infty} \right)\|\phi\|_{Y^\ast}.
		\]
		Finally, we consider $Z = H_Y$. By using the definition of $H_Y$, Definition~\Rref{defi:weighted_sobolev}, and then applying a change of variables \cite[Theorem 3.9]{EvaG92}, we conclude
		\begin{align*}
			& \| S^\ast(t) \phi \|^2_{H_Y} - \| S^\ast(t) \phi \|^2_{H_X} = \int_\Omega \nabla \phi(\Phi^t(z))^\top \mathrm{D} \Phi^t(z)^\top \mathrm{D} \Phi^t(z) \nabla \phi ( \Phi^t(z)) \dint z \\
			\le & \int_{\Phi^t(\Omega)} |\det{ \mathrm{D} \Phi^t( \Phi^{-t}(a)) } |^{-1} \nabla \phi(a)^\top \mathrm{D} \Phi^t(\Phi^{-t}(a))^\top \mathrm{D} \Phi^t( \Phi^{-t}(a)) \nabla \phi ( a) \dint a.
		\end{align*}
		Applying \eqref{eq:prop:bound_semigroup:boundDphi} then leads to
		\[
			\| S^\ast(t) \phi \|^2_{H_Y} \le \| S^\ast(t) \phi \|^2_{H_X} + C_2^2 L^{d} \| \nabla \phi \|_{\Lp{2}{\Omega}}^2 \le L^{d}\left( \exp(2 \omega_0 t) + C_2^2 \right) \| \phi\|_{\Wpw{2}{w}{\Omega}}^2,
		\]
		where we used the result from the case $Z = H_X$ for the last inequality. Therefore,
		\[
		\| S^\ast(t)\|_{\mathcal{L}(H_Y)} \le L^{\sfrac{d}{2}}\sqrt{\exp(2\omega_0 t)+C_2^2},
		\]
		which was to be shown.
	\end{proof}
	\begin{thrm}
		\label{thm:cont_semigroup}
		For $Z \in \{ X, H_X, Y, H_Y \}$ and a boundary-admissible dynamic $f$, the Koopman semigroup associated with $f$ is a well-defined weak-* continuous semigroup, while the Perron--Frobenius semigroup associated with $f$ defines a strongly continuous semigroup.
	\end{thrm}
	\begin{proof}
		Since $f$ is autonomous, it follows that $\Phi^t(\Phi^s(z)) = \Phi^{t+s}(z)$ and, consequently, $S^\ast(t+s) = S^\ast(t)S^\ast(s)$.
		It remains to show weak-* continuity of $S^\ast(t)$. We begin with $Z := X = \fLpw{1}{w}{\Omega}$. Let $\phi \in \fWpw{\infty}{w}{\Omega}$. Then, with Definition~\Rref{defi:preadjoint_semigroup}, we can rewrite
		\[
		(S(t) \phi - \phi)(x) = \left\{ 
		\begin{array}{lcl}
			\mu_t(x) \phi( \Phi^{-t} (x) )  - \phi(x) & \qquad &  \text{for} \; x \in \Phi^t(\Omega) \\
			-\phi(x) & \qquad & \text{else}
		\end{array}
		\right. .
		\]
		From Proposition~\Rref{prop:bounds_flow} and Lemma~\Rref{lemma:dissipation_type}, we obtain the bound
		\[
		\| S(t) \phi - \phi \|_{\Lpw{\infty}{w}{\Omega}} \le ( 1 + \| \mu_t \|_{\Lp{\infty}{\Omega}} ) \| \phi \|_{\Lpw{\infty}{w}{\Omega}} \le ( 1 + \exp(2 \omega_0 \delta ) L_\delta^{d}  ) \| \phi \|_{\Lpw{\infty}{w}{\Omega}},
		\]
		for $0 \le t < \delta$, where $L_\delta:=\exp(\|f\|_{\Wpw{\infty}{w}{\Omega}}\delta)$. 
		Since $\Omega$ is bounded, any sequence $\{ (S(t_n) \phi - \phi) w \}_{ n \in \N}$ with $0 \le t_n \le \delta$ is therefore bounded by an integrable function almost everywhere. We can apply the dominated convergence theorem \cite[Theorem A3.21]{Alt12} to conclude
		\begin{align*}
			& \lim_{t \rightarrow 0} \int_\Omega | (S(t) \phi)(x) - \phi(x) | w(x) \dint x \\
			= & \int_\Omega \lim_{t \rightarrow 0} \left| \chi_{\Phi^t(\Omega)} \frac{\phi(\Phi^{-t}(x)) w^2(\Phi^{-t}(x))}{|\det{\mathrm{D} \Phi^{t}(\Phi^{-t}(x))}| w^2(x)} - \phi(x) \right| w(x) \dint x  = 0,
		\end{align*}
		where we used the continuity of $\phi$, Proposition~\Rref{prop:bounds_flow}, and Lemma~\Rref{lemma:dissipation_type}. We conclude that
		\begin{align*}
			\lim_{t \rightarrow 0} \| S(t) \phi - \phi\|_{X} &= 0, \\
			\lim_{t \rightarrow 0} \| S(t) \phi - \phi\|_{H_X} &= 0,
		\end{align*}
		where the second limit follows by applying the same dominated convergence argument to the squared integrands.
		Therefore, for $\psi \in X^\ast$, we have shown that
		\[
		\lim_{t \rightarrow 0} |\inner[X,X^\ast]{\phi}{S^\ast(t) \psi - \psi }| \le \lim_{t \rightarrow 0} \| S(t)\phi - \phi \|_X \| \psi \|_{X^\ast} =  0 \qquad \text{for all}\; \phi \in D, 
		\]
		where $D = \fWpw{\infty}{w}{\Omega} \subseteq X$ is dense in $X$ by Lemma~\Rref{lemma:density_wpwinfty}. By Lemma~\Rref{lemma:weak_densesubset}, it follows that $S^\ast$ is a weak-* continuous semigroup on $X^\ast$ and similarly on $H_X$.\\
		Next, consider $[(g,0,\dots,0)]_\sim \in D$ and $\phi \in Y^\ast$, where $D$ is defined as in Lemma~\Rref{lemma:dense_zeros} and is dense in $Y$. Then
		\[
		\lim_{t \rightarrow 0} | \inner[Y,Y^\ast]{[(g,0,\dots,0)]_\sim}{S^\ast(t)\phi - \phi } | = \lim_{t \rightarrow 0} | \inner[X,X^\ast]{g}{S^\ast(t) \phi - \phi} | = 0.
		\]
		Again, with Lemma~\Rref{lemma:weak_densesubset}, we conclude that $S^\ast$ is a weak-* continuous semigroup for $Z=Y$. The case $Z = H_Y=\fWpw{2}{w}{\Omega}$ follows similarly.\\
		Next, we show that $S(t)$ is a strongly continuous semigroup. Let $\phi \in Z$ and $\psi \in Z^\ast$. Then, by the weak-* continuity of $S^\ast(t)$,
		\[
		\lim_{t \searrow 0} \inner[Z,Z^\ast]{S(t)\phi}{\psi} = 	\lim_{t \searrow 0} \inner[Z,Z^\ast]{\phi}{S^\ast(t)\psi} = \inner[Z,Z^\ast]{\phi}{\psi},
		\]
		which yields that $S(t)$ is weakly continuous. By \cite[Theorem 5.8]{EngN06}, $S(t)$ is also strongly continuous.
	\end{proof}
	\begin{prpstn}
		\label{prop:generator_koopman}
		Let $f$ be a boundary-admissible dynamic and let $S^\ast$ be the corresponding Koopman semigroup.
		For $Z\in\{X,H_X, Y, H_Y\}$, denote by $A_Z^\ast$ the weak-$\ast$
		generator of $S^\ast$ on $Z^\ast$ with domain
		\[
		\mathcal D(A_Z^\ast) := \left\{
			\phi\in Z^\ast \; \big| \; \wslim_{t\searrow0} \tfrac{1}{t} \left( S^\ast(t) \phi - \phi \right) \text{ exists in } Z^\ast
		\right\}.
		\]
		Then
		\[
			A_Z^\ast\phi=f^\top\nabla\phi
			\qquad\text{with}\qquad
		\mathcal D(A_Z^\ast)
		=
		\left\{
			\phi\in Z^\ast \; \big|\; f^\top \nabla \phi\in Z^\ast
		\right\}.
		\]
		Here $f^\top\nabla\phi=g\in Z^\ast$ is understood in a distributional sense, i.e.,
		\begin{equation}
			-\int_\Omega \phi\,\operatorname{div}(\psi f)\dint x
			=
			\int_\Omega g\psi \dint x
			\qquad\text{for all }\psi\in \CcInfty{\Omega}.
			\label{eq:weak_fTnabla}
		\end{equation}
	\end{prpstn}
	\begin{proof}
		First, consider $Z\in\{X,H_X\}$. Let $\phi\in\mathcal D(A_Z^\ast)$. Then there exists $g\in Z^\ast$ such that
		\[
			\lim_{t\searrow0}
			\inner[Z,Z^\ast]{\rho}{\tfrac{1}{t}\left(S^\ast(t)\phi-\phi\right)}
			=\inner[Z,Z^\ast]{\rho}{g}
			\qquad\text{for all }\rho\in Z.
		\]
		In particular, for every $\psi\in\CcInfty{\Omega}$,
		\[
			\lim_{t\searrow0}
			\inner[Z,Z^\ast]{\psi w^{-1}}{\tfrac{1}{t}\left(S^\ast(t)\phi-\phi\right)}
			=\inner[Z,Z^\ast]{\psi w^{-1}}{g}.
		\]
		By Lemma~\Rref{lemma:preadjoint_semigroup} and Definition~\Rref{defi:preadjoint_semigroup},
		\begin{align*}
			&\lim_{t\searrow0}
			\inner[Z,Z^\ast]{\psi w^{-1}}{\tfrac{1}{t}\left(S^\ast(t)\phi-\phi\right)}
			 =\lim_{t\searrow0}
			\inner[Z,Z^\ast]{\tfrac{1}{t}\left(S(t)(\psi w^{-1})-\psi w^{-1}\right)}{\phi}\\
			= &\lim_{t\searrow0}\tfrac{1}{t}\int_\Omega\phi(x)
			\left(
				\chi_{\Phi^t(\Omega)}(x)
				\frac{\psi(\Phi^{-t}(x))w(\Phi^{-t}(x))}
				{\left|\det{\mathrm D\Phi^t(\Phi^{-t}(x))}\right|}
				-\psi(x)w(x)
			\right)\dint x.
		\end{align*}
		Proposition~\Rref{prop:explicit_Dphi}, Liouville's formula~\cite[Exercise 2.2]{AmbC14}, and the chain rule~\cite[Theorem 3.3]{EvaG92} yield, for almost every $a\in\Phi^t(\Omega)$,
		\[
			\frac{\mathrm d}{\mathrm dt}
			\left|\det{\mathrm D\Phi^t(\Phi^{-t}(a))}\right|^{-1}
			=
			-\operatorname{div}(f)(\Phi^{-t}(a))
			\left|\det{\mathrm D\Phi^t(\Phi^{-t}(a))}\right|^{-1}
		\]
		and
		\[
			\frac{\mathrm d}{\mathrm dt}
			\frac{w(\Phi^{-t}(a))}{w(a)}
			=
			\frac{w(\Phi^{-t}(a))}{w(a)}
			(fw)^\top\nabla(w^{-1})(\Phi^{-t}(a)).
		\]
		Since $\psi$ has compact support, there exists $h>0$ such that
		$\operatorname{supp}\psi\subseteq\Phi^t(\Omega)$ for all $0\le t\le h$.
		Using the preceding identities, the Fubini--Tonelli theorem~\cite[Remark A6.11]{Alt12}, and the area formula~\cite[Theorem 3.8]{EvaG92}, we obtain
		\begin{align*}
			&\lim_{t\searrow0}
			\inner[Z,Z^\ast]{\psi w^{-1}}{\tfrac{1}{t}\left(S^\ast(t)\phi-\phi\right)}
			=
			-\lim_{t\searrow0}\tfrac{1}{t}\int_0^t
			\int_{\Phi^s(\Omega)}
			\phi(x)
			\frac{\operatorname{div}(fw\psi)(\Phi^{-s}(x))}
			{\left|\det{\mathrm D\Phi^s(\Phi^{-s}(x))}\right|}
			\dint x\dint s\\
			= & -\lim_{t\searrow0}\tfrac{1}{t}\int_0^t
			\int_\Omega \phi(\Phi^s(a))\operatorname{div}(fw\psi)(a)
			\dint a\dint s\\
			= & -\lim_{t\searrow0} \tfrac{1}{t} \int_0^t \inner[Z,Z^\ast]{w^{-2}\operatorname{div}(fw\psi)}{S^\ast(s)\phi} \dint s 
			= - \inner[Z,Z^\ast]{w^{-2}\operatorname{div}(fw\psi)}{\phi}  = -\int_\Omega\phi\,\operatorname{div}(fw\psi)\dint x.
		\end{align*}
		Here, note that $w^{-2}\operatorname{div}(fw\psi)\in Z$ and that the convergence as $t\searrow0$ follows from Theorem~\Rref{thm:cont_semigroup}. Comparing the limit with the definition of $g$, we obtain
		\[
			-\int_\Omega \phi\,\operatorname{div}(\psi f w)\dint x
			=
			\int_\Omega g\psi w\dint x
			\qquad\text{for all }\psi\in\CcInfty{\Omega}.
		\]
		For $p=1$ when $Z=X$ and $p=2$ when $Z=H_X$, note that
		\[
			\frac{\operatorname{div}(fw\psi)}{w}
			=
			\left(\operatorname{div}f-(fw)^\top\nabla(w^{-1})\right)\psi
			+f^\top\nabla\psi.
		\]
		The coefficients on the right-hand side are bounded. Hence both sides of the equality are continuous with respect to the $W^{1,p}(\Omega)$ norm of $\psi$, and the identity extends to every $\psi\in W_0^{1,p}(\Omega)$. For any $\eta\in\CcInfty{\Omega}$, we have $w^{-1}\eta\in W_0^{1,p}(\Omega)$ because $w^{-1}\in W^{1,\infty}(\Omega)$. Taking $\psi=w^{-1}\eta$ gives
		\[
			-\int_\Omega\phi\,\operatorname{div}(\eta f)\dint x
			=
			\int_\Omega g\eta\dint x,
		\]
		which is \eqref{eq:weak_fTnabla}. For $Z\in\{Y,H_Y\}$, one uses the dense sets of zeroth-component tests
		\[
			D_Z:=\left\{[(\psi w^{-1},0,\dots,0)]_\sim
			\mid\psi\in\CcInfty{\Omega}\right\}\subset Z.
		\]
		Their density follows from Lemma~\Rref{lemma:dense_zeros} and the preceding density argument for $X$ and $H_X$. The same calculation applies to these tests. In the limit above, one uses $Y^\ast\subset X^\ast$ and $H_Y\subset H_X$, together with Theorem~\Rref{thm:cont_semigroup} on $X^\ast$ and $H_X$, respectively. Thus the forward implication also holds for $Y$ and $H_Y$.

		Conversely, let $Z\in\{X,H_X,Y,H_Y\}$, and suppose that $\phi\in Z^\ast$ and $g:=f^\top\nabla\phi\in Z^\ast$ satisfy \eqref{eq:weak_fTnabla}. Since $f$ is autonomous and boundary-admissible, its trajectories remain in $\Omega$, i.e., $\Phi^s(x)\in\Omega$ for all $x\in\Omega$ and $s\ge0$. Therefore, the argument of Proposition~2.3 of \cite{AmbC14}, localized to $\Omega$, yields
		\[
			S^\ast(t)\phi-\phi=\int_0^t S^\ast(s)g\dint s
			\qquad\text{almost everywhere in }\Omega,
		\]
		where the integral is understood in the weak-$\ast$ sense. For any fixed $t_0>0$, Proposition~\Rref{prop:bound_semigroup} therefore gives
		\[
			\sup_{0<t\le t_0}
			\left\|\tfrac{1}{t} \left( S^\ast(t)\phi-\phi \right) \right\|_{Z^\ast}
			\le
			\sup_{0\le s\le t_0}\|S^\ast(s)\|_{\mathcal L(Z^\ast)}
			\|g\|_{Z^\ast}
			=:M_{t_0}\|g\|_{Z^\ast}.
		\]
		Consider the dense sets
		\[
			D_Z:=
			\begin{cases}
				\{\psi w^{-1}\mid\psi\in\CcInfty{\Omega}\},
				& Z\in\{X,H_X\},\\
				\{[(\psi w^{-1},0,\dots,0)]_\sim
				\mid\psi\in\CcInfty{\Omega}\},
				& Z\in\{Y,H_Y\}.
			\end{cases}
		\]
		For every $\rho_0\in D_Z$, the flow identity and Theorem~\Rref{thm:cont_semigroup} yield
		\[
			\lim_{t\searrow0}
			\inner[Z,Z^\ast]{\rho_0}
			{\tfrac{1}{t} \left( S^\ast(t) \phi - \phi \right)-g}
			=
			\lim_{t\searrow0}\frac{1}{t}\int_0^t
			\inner[Z,Z^\ast]{\rho_0}{S^\ast(s)g-g}\dint s
			=0.
		\]
		Now let $\rho\in Z$. For every $\varepsilon>0$, choose $\rho_0\in D_Z$ such that $\|\rho-\rho_0\|_Z<\varepsilon$. Then
		\begin{align*}
			&\limsup_{t\searrow0}
			\left|
			\inner[Z,Z^\ast]{\rho}
			{\tfrac{1}{t} \left( S^\ast(t) \phi - \phi \right)-g}
			\right|
			\le
			(M_{t_0}+1)\|g\|_{Z^\ast}\|\rho-\rho_0\|_Z
			\le
			(M_{t_0}+1)\|g\|_{Z^\ast}\varepsilon.
		\end{align*}
		Letting $\varepsilon\searrow0$, we conclude that
		\[
			\wslim_{t\searrow0}
			\tfrac{1}{t} \left( S^\ast(t) \phi - \phi \right)
			=g=f^\top\nabla\phi.
		\]
		Hence $\phi\in\mathcal D(A_Z^\ast)$ and $A_Z^\ast\phi=f^\top\nabla\phi$.
	\end{proof}
	\subsection{Stability of the Koopman semigroup}
	\label{ssec:stability_koopman}
	 To investigate stability, we assume that $0\in\Omega$ and set $w(x) = \sfrac{1}{\|x\|}$. It turns out that this weighting is suitable for systems whose linearization at the origin is stabilizable. While stability of the Koopman semigroup is investigated in \cite{MauSM20stab}, our analysis is carried out over a different space and hence with respect to a different norm. We start with the definition of a compatible dynamic.
	 \begin{dfntn}
	 	\label{defi:linear_stable}
		Let $w = \| \cdot\|^{-1}$. A function $f \colon \Omega \to \R^d$ is called a \emph{linear stable dynamic} if the following conditions are satisfied:
	 	\begin{enumerate}[label=(\roman*)]
	 		\item $f$ is a boundary-admissible dynamic, see Definition~\Rref{defi:bd_admissible_dynamic}.
	 		\label{defi:stability:dynamic}
			\item $\mathrm{D}f$ is continuous at $0$ and there exists a positive definite solution $\hat{P}$ to the Lyapunov equation
	 		\[
	 		\mathrm{D}f(0)^\top \hat{P} + \hat{P} \mathrm{D}f(0) + \hat{H} = 0
	 		\]
	 		for some positive definite $\hat{H} \in \R^{d\times d}$.
	 		\label{defi:stability:lqr}
	 	\end{enumerate}
	 \end{dfntn}
	The Koopman semigroup associated with a linear stable dynamic exhibits exponential stability when restricted to a sufficiently small neighborhood around the equilibrium point.
	\begin{prpstn}
		\label{prop:local_uncontrolled_stability}
		Let $w = \| \cdot \|^{-1}$, let $f$ be a linear stable dynamic, and let $\hat{H}$ and $\hat{P}$ denote the positive definite matrices from Definition~\Rref{defi:linear_stable} \Rref{defi:stability:lqr}. For
		\[
		\Omega_\delta := \{ z \in \Omega \; \text{such that}\; \| z \|_{\hat{P}} := (z^\top \hat{P} z)^{\sfrac{1}{2}} < \delta \},
		\]
		define $Y_\delta :=\fNWpw{1}{w}{ \Omega_\delta }$. Then the following hold:
		\begin{enumerate}[label=(\roman*)]
			\item For all sufficiently small $\delta > 0$, the restriction $S_\delta(t)\colon Y_\delta \to Y_\delta$ is well-defined.
			\item For all $\omega_1 > \omega_0$ there exists $\delta > 0$ such that 
			\label{prop:local_uncontrolled_stability:stability}
			\[
			\| S_\delta(t) \|_{\mathcal{L}(Y_\delta)} \le C(\omega_1) \exp( t \omega_1 ) \qquad \text{for all}\; t \in (0,\infty),
			\]
			where $\omega_0 = - \frac{1}{2} \min_{1 \le i \le d} \lambda_i( \hat{P}^{-\sfrac{1}{2}} \hat{H} \hat{P}^{-\sfrac{1}{2}} )$.
		\end{enumerate}
	\end{prpstn}
	\begin{proof}	
		First, we show that the semigroup is well-defined over $Y_{\delta}$ for $\delta >0$ small enough by showing that $f\big|_{\Omega_\delta}$ is a boundary-admissible dynamic.
		 The first property from Definition~\Rref{defi:bd_admissible_dynamic}\Rref{defi:bd_admissible_dynamic_regularity} follows directly from $f \in \fWpw{\infty}{w}{\Omega}$. For Definition~\Rref{defi:bd_admissible_dynamic}\Rref{defi:bd_admissible_dynamic_tangent}, note that the outward unit normal to $\Omega_\delta$ is $\nu(z) := \tfrac{\hat{P} z}{\| \hat{P} z\|}$ for $z \in \partial \Omega_\delta$. Since $f$ is differentiable at $0$ and $f(0) = 0$, which is a consequence of the weighting, we conclude
		\begin{align*}
			&	\| \hat{P} z \| \nu(z)^\top f(z) = z^\top \hat{P} f(z) \\
			= & \tfrac{1}{2}z^\top \left(  \hat{P}\mathrm{D} f( 0 ) + \mathrm{D} f(0)^\top \hat{P} + o(1) \right) z = - \tfrac{1}{2} z^\top( \hat{H} +
			o(1)) z.
		\end{align*}
		Since this is negative for small enough $\delta$, we conclude that both properties are satisfied and $f\big|_{\Omega_\delta}$ is a boundary-admissible dynamic, which allows us to apply Theorem~\Rref{thm:cont_semigroup} to the restricted system.
		
		To show \Rref{prop:local_uncontrolled_stability:stability}, let $\omega_1 > \omega_0$, $\phi \in Y_\delta$, and $\varepsilon >0$. By definition, there exist $\phi_0,\phi_1, \dots, \phi_d$ such that
		\[
		(\phi_0,\phi_1, \dots, \phi_d) \in [\phi]_{\sim} \quad \text{and} \quad  \| \phi_0 \|_{\Lpw{1}{w}{\Omega_\delta}} + \sum_{i=1}^d \| \phi_i \|_{\Lp{1}{\Omega_\delta}} \le \| \phi \|_{Y_\delta}  + \varepsilon.
		\]
		By the characterization of the semigroup in Lemma~\Rref{lemma:preadjoint_semigroup}, we conclude
		\begin{align*}
			\| S_\delta(t) \phi \|_{Y_\delta} \le & \| S_\delta(t) \phi_0 \|_{\Lpw{1}{w}{\Omega_\delta}} + \int_{\Phi^t(\Omega_\delta)} \| M_t(a) \left( \phi_1(\Phi^{-t}(a)), \dots, \phi_d( \Phi^{-t}(a)) \right)^\top \|_1 \dint a.
		\end{align*}
		For the first part, we can use the same argument given in \cite[Theorem 2.8]{BreH23} to conclude
		\[
		\| S_\delta(t) \phi_0 \|_{\Lpw{1}{w}{\Omega_\delta}} \le C  \exp( \omega_1 t ) \| \phi_0 \|_{\Lpw{1}{w}{\Omega_\delta}}.
		\]
		Next, for some $p \in \R^d$ with $\|p\| = 1$ and for almost every $z \in \Omega_\delta$ we consider
		\[
		r_{z,p} \colon [0,\infty) \to \R_+, \quad t \mapsto p^\top \mathrm{D}\Phi^t(z)^\top \hat{P} \mathrm{D} \Phi^t(z) p.
		\]
		By Proposition~\Rref{prop:explicit_Dphi}, $r_{z,p}$ is the solution to
		\[
		\left\{
		\begin{array}{rcl}
			\tfrac{\mathrm{d}}{\mathrm{d} t} r_{z,p}(t) &=& p^\top \,\mathrm{D}\Phi^t(z)^\top \left(  \mathrm{D}f(\Phi^t(z))^\top \hat{P}  + \hat{P} \mathrm{D}f( \Phi^t(z) ) \right) \mathrm{D}\Phi^t(z) p\\
			r_{z,p}(0) & =& p^\top \hat{P} p
		\end{array}
		\right.
		\]
		in the sense of Carath\'{e}odory. Since $\mathrm{D}f$ is continuous at $0$ and $\hat{P}$ is the solution to the algebraic Lyapunov equation, we have
		\[
		\mathrm{D}f(\Phi^t(z))^\top \hat{P} + \hat{P} \mathrm{D}f(\Phi^t(z) ) 
		\preceq   - ( 1 - o(1)) \hat{H}
		\preceq 2 ( 1 - o(1)) \omega_0 \hat{P},
		\]
		where $o(1)\to0$ as $\delta\searrow0$ and we used $-\hat{H} \preceq 2 \omega_0 \hat{P}$ by the definition of $\omega_0$.
		Therefore, $\frac{\mathrm{d}}{\mathrm{d} t} r_{z,p}(t) \le 2 ( 1- o(1) ) \omega_0 r_{z,p}(t)$. By choosing $\delta >0$ small enough, we conclude from Gronwall's lemma~\cite[§29, Theorem VI]{Wal98} that $r_{z,p}(t) \le r_{z,p}(0)\exp( 2 \omega_1 t )$ for any $\omega_1 > \omega_0$. Since $\hat{P}$ is invertible by assumption, we conclude for almost every $z \in \Omega_\delta$
		\[
		\| \mathrm{D} \Phi^t(z) \|_2 \le C \sup_{\|p\|_2=1} \| \mathrm{D} \Phi^t(z) p \|_{\hat{P}} \le C \exp( \omega_1 t ).
		\]
		Finally, we use the explicit expression for $M_t$ from Lemma~\Rref{lemma:preadjoint_semigroup} to compute
		\begin{align*}
			&\int_{\Phi^t(\Omega_\delta)} \| M_t(a) \left( \phi_1(\Phi^{-t}(a)), \dots, \phi_d( \Phi^{-t}(a)) \right)^\top \|_1 \dint a \\
			= &\int_{\Omega_\delta} \| \mathrm{D}\Phi^t(x) \left( \phi_1(x), \dots, \phi_d(x) \right)^\top \|_1 \dint x\\
			\le & C \exp(\omega_1 t ) \int_{\Omega_\delta} \| \left( \phi_1(x), \dots, \phi_d(x) \right)^\top \|_1 \dint x
			\le  C \exp(\omega_1 t ) \sum_{i=1}^d \| \phi_i \|_{\Lp{1}{\Omega_\delta}}.
		\end{align*}
		Since the choice of $\delta$ was independent of $\phi$ and $\varepsilon$, we conclude that for every $\varepsilon > 0$ and fixed $t > 0$,
		\[
		\| S_\delta(t) \phi \|_{Y_\delta} \le C \exp( \omega_1 t ) \left( \| \phi \|_{Y_\delta} + \varepsilon \right).
		\]
		Letting $\varepsilon \to 0$ proves the estimate.
	\end{proof}
	\section{The value function as a bilinear form}
	\label{ssec:koopman_feedback}
	We now extend the theory to controlled dynamical systems with a single input. 
	We do so by fixing a feedback law and analyzing the closed-loop dynamics with the tools from Section~\Rref{sec:koopman}. Although the generalization to multiple inputs is straightforward, we do not pursue it here in order to keep the analysis as simple as possible. We begin with the basic definitions.
	\begin{dfntn}
		Let $w := \| \cdot \|^{-1}$. 
		\begin{enumerate}[label=(\roman*)]
			\item A tuple $(f,b)$ is called a \emph{boundary-admissible controlled dynamic pair} if $f$ is a boundary-admissible dynamic and $b_i \in \fWp{\infty}{\Omega}$ for $1\le i\le d$. 
			\item A sequence $\{ c_i \}_{i\in \N} \subset \fWpw{\infty}{w}{\Omega}$ with $\sum_{i=1}^\infty \| c_i \|^2_{Y^\ast}<\infty$ is called a \emph{nuclear sequence of observables}. 
			\item A state penalty $g \colon \Omega \to \R_+$ is called a \emph{nuclear cost} if there exists a nuclear sequence of observables $\{c_i\}_{i \in \N}$ such that
			\[
				g(z) = \sum_{i=1}^\infty c_i(z)^2 \qquad \text{for all} \; z \in \Omega
			\]
			and a \emph{finite-rank cost} if $c_i \equiv 0$ for all $i > r$ and some $r \in \N$.
		\end{enumerate}
	\end{dfntn}
	For such systems we can define a value function for a state-constrained control problem, which requires the controlled trajectories to stay inside $\Omega$; see, e.g., \cite{CapL90}. This is equivalent to imposing a boundary condition for the HJB equation, and its solution can be used to derive an optimal feedback. The following proposition is a direct consequence of the Carath\'{e}odory existence and uniqueness theorem \cite[§10, Theorem VI]{Wal98}. 
	\begin{prpstn}
		\label{prop:picard_lindeloeff}
		For a given $\alpha \in \mathrm{L}^\infty_{\mathrm{loc}}(0,\infty;\R)$ and $z \in \Omega$, there exists a maximal time horizon $T_z(\alpha) \in [0,\infty) \cup \{ \infty\}$ such that
		\eqref{eq:dynamics_intro} has a unique Carath\'{e}odory solution $x(t;\alpha) \in \overline{\Omega}$ on $[0,T_z(\alpha))$ with $x(0;\alpha)=z$.
	\end{prpstn}
	\begin{dfntn}
		\label{defi:value_function}
		Let $(f,b)$ be a boundary-admissible controlled dynamic pair and $\{c_i\}_{i \in \N}$ a nuclear sequence of observables. We define the value function with state constraints as
		\[
		v(z) := \inf_{\alpha \in \mathcal{A}_{z} } \; \int_0^\infty \left( \sum_{i=1}^\infty c_i( x( t ;\alpha) )^2 + \alpha(t)^2 \right)
		\dint{t}
		\]
		where $x(t;\alpha)$ denotes the solution of the ODE in \eqref{eq:dynamics_intro} and the set of admissible control inputs is defined as
		\[
			\mathcal{A}_z := \{ \alpha \in \mathrm{L}^\infty_{\mathrm{loc}}(0,\infty;\R) \; \text{such that}\; T_z(\alpha) = \infty \} \ni 0.
		\]
	\end{dfntn}
	The time-dependent formulation does not allow us to apply the results from the previous sections directly. Therefore, we introduce feedback laws $u \in\fWpw{\infty}{w}{\Omega}$ with associated closed-loop trajectories:
	\begin{align}
		\label{eq:feedback_controlled}
		\left\{
		\begin{array}{rclrl}
			\tfrac{\mathrm{d}}{\mathrm{d}t} x_u(t) &=& f(x_u(t)) + b(x_u(t)) u(x_u(t)) & \qquad \text{for all}\; & t \in (0,\infty) \\
			x_u(0) &=& z,& &
		\end{array}\right.
	\end{align}
	for some $(f,b)$ and initial value $z \in \Omega$. In analogy to Definition~\Rref{defi:uncontrolled_flow}, we define the \emph{feedback flow} as $\Phi_u^t(z) := x_u(t)$.
	As in the uncontrolled case, trajectories must remain in $\Omega$; otherwise, the closed-loop Koopman semigroup is not well-defined. To enforce this property for the closed-loop trajectories, we restrict the set of feedback laws under consideration as follows.
	\begin{dfntn}
		\label{defi:boundary_compatible}
		Let $(f,b)$ be a boundary-admissible controlled dynamic pair. A feedback $u \in Y^\ast$ is called \emph{boundary-admissible} if 
		\[
		\nu(x)^\top \left( f(x) + b(x) u(x) \right) \le 0 \qquad \text{for all}\; x \in \partial \Omega.
		\]
		The triple $(f,b,u)$ is called a \emph{boundary-admissible closed-loop dynamic}.
	\end{dfntn}
	\begin{rmrk}
		\label{rem:vanishing_boundary}
		If $b$ vanishes on the boundary, i.e., $b\big|_{\partial\Omega} \equiv 0$, then any feedback $u \in Y^\ast$ is automatically boundary-admissible. 
		This can be interpreted as imposing an increasingly strong control penalty near the boundary, akin to a barrier method. 
	\end{rmrk}
	As is well known, optimal control allows for feedback formulations by means of the gradient of the value function under the assumption that the value function is sufficiently regular.
	\begin{thrm}
		\label{thm:optimal_feedback}
		Let $(f,b)$ be a boundary-admissible controlled dynamic pair and $\{c_i\}_{i \in \N}$ a nuclear sequence of observables. If $v \in W^{2,\infty}(\Omega)$, then $v$ solves the HJB equation in the strong sense:
		\begin{align*}
			0 = \min_{\beta \in \R} \left( \sum_{i=1}^\infty c_i(z)^2 + \beta^2 + (f(z) + \beta b(z))^\top \nabla v(z) \right)
		\end{align*}
		for every $z \in \Omega$.
		Furthermore, the feedback control
		\[
		u_* \colon  \Omega  \to \mathbb{R}, \quad  z \mapsto - \tfrac{1}{2} b(z)^\top \nabla v(z)
		\]
		satisfies $(f(z) + u_*(z) b(z) )^\top \nu(z) \le 0$ for all $z \in \partial \Omega$ and is optimal in the sense that for any $z \in \Omega$, the control attains the infimum, i.e.,
		\[
		v(z) = \int_0^\infty \left( \sum_{i=1}^\infty c_i( \Phi_{u_*}^t(z) )^2 + \alpha(t)^2 \right) \, \dint{t} \qquad \text{with} \quad  \alpha(t) := u_*( \Phi_{u_*}^t(z) ).
		\]
	\end{thrm}
	\begin{proof}
		 See Appendix~\Rref{sec:technical_proofs}.
	\end{proof}
	\begin{rmrk}
		The infinite-horizon problem without a discount factor is often not covered in the HJB literature because the HJB equation does not necessarily admit a viscosity solution.
		For discounted costs with the same state constraints, a similar result, together with uniqueness and existence results, can be found in \cite{CapL90,CraL86}. Note that discounting would translate into a shifted semigroup whose generator differs by a scalar multiple of the identity.
		Assuming $v \in W^{2,\infty}(\Omega)$ is relatively restrictive and leads to strong solutions rather than viscosity solutions of the HJB equation. However, the assumption is necessary in our case since the derived feedback $u_*$ needs to be an element of $\fWpw{\infty}{w}{\Omega} = Y^\ast$. 
	\end{rmrk}
	
	Our goal now is to show that the value function can be linked to a bilinear form that allows for finite-rank approximations, leading to a sum-of-squares approximation of the value function. We thus introduce the Koopman semigroup associated with the closed-loop dynamics in order to analyze how the feedback and the dynamics interact.
	\begin{dfntn}
		\label{defi:composition_semigroup:control}
		Let $(f,b,u)$ be a boundary-admissible closed-loop dynamic. We define the \emph{closed-loop Koopman semigroup} as 
		\[
		S_u^\ast \colon [0,\infty)  \to \mathcal{L}(Y^\ast), \quad t \mapsto S^\ast_u(t)
		\]
		where $S_u^\ast(t) \phi := \phi \circ \Phi_{u}^t$ for all $ \phi \in Y^\ast$.
	\end{dfntn}
	\begin{rmrk}
		For a fixed, boundary-admissible feedback $u \in Y^\ast$, we can apply the results from Sections~\Rref{ssec:koopman_uncontrolled} and~\Rref{ssec:stability_koopman} since the closed-loop dynamic $f + bu$ is boundary-admissible. In particular, Theorem~\Rref{thm:cont_semigroup} shows that the controlled composition semigroup is indeed a weak-* continuous semigroup. Furthermore, this yields the existence of a preadjoint $S_u$ and its explicit representation from Lemma~\Rref{lemma:preadjoint_semigroup}.
	\end{rmrk}
	Let us incorporate the control penalty into a suitable observation operator.
	\begin{dfntn}
		\label{defi:control_operator}
		Let $(f,b,u)$ be a boundary-admissible closed-loop dynamic and $\{c_i\}_{i\in \N} \subset \fWpw{\infty}{w}{\Omega}$ a nuclear sequence of observables. We define the \emph{observation with control penalty} $C_u$ as the operator
		\[
		C_u \colon Y  \to \ell_2, \quad \phi \mapsto ( \inner[Y,Y^\ast]{\phi}{u} , \inner[Y,Y^\ast]{\phi}{c_1}, \dots ).
		\]
		To simplify notation we will set $c_0 := u$.
	\end{dfntn}
	The following result extends \cite[Proposition 4.2]{BreH23} to the controlled case.
	\begin{lmm}
		\label{lemma:exponentially_stability}
		Let $(f,b)$ be a boundary-admissible controlled dynamic pair, and suppose that the following assumptions are satisfied:
		\begin{enumerate}[label=(\roman*)]
			\item The linearized pair $(\mathrm{D}f(0),b(0))$ is stabilizable.
			\item There are finitely many observables $c_1,\dots,c_r$ and 
			\[
			\label{req:exp_stab:rank}
			\hat{C}:=(\nabla c_1(0), \dots,\nabla c_r(0))^\top,
			\qquad \operatorname{rank}(\hat{C})=d.
			\]
			\item For every $\delta > 0$, there exists $\varepsilon > 0$ such that the collective observables satisfy
			\[
			\label{req:exp_stab:positivity}
			\sum_{i=1}^r c_i(z)^2 \ge \varepsilon \qquad \text{for all} \; z \in \Omega \setminus B_\delta(0).
			\]
			\item The value function satisfies $v \in \mathrm{W}^{2,\infty}(\Omega)$ and $\mathrm{D}^2 v$ is continuous at $0$.
			\label{req:exp_stab:value_function}
		\end{enumerate}
		Then $u_*$ is a boundary-admissible feedback, and the preadjoint $S_{u_*}(t)$ of the semigroup from Definition~\Rref{defi:composition_semigroup:control} is exponentially stable over $Z \in \{ X, Y\}$ with type at most
		\begin{align*}
			\omega_0 := & \inf \; \{ \omega \in \R \;  \vert \; \sup_{ t\ge 0} \{ \exp( - \omega t) \| S_{u_*}(t)\|_{\mathcal L(Z)} \} < \infty \} \\
			\le  &\tfrac{1}{2}\max \{ - \lambda_i(\hat{P}^{-\sfrac{1}{2}} \hat{C}^\top \hat{C} \hat{P}^{-\sfrac{1}{2}}) \; \vert \; 1 \le i \le d\} < 0,
		\end{align*}
		where $\hat{P}$ solves the algebraic Riccati equation
		\[
		\mathrm{D} f(0)^\top \hat{P} + \hat{P} \mathrm{D} f(0) - \hat{P} b(0) b(0)^\top \hat{P} + \hat{C}^\top \hat{C} = 0.
		\]
	\end{lmm}
	\begin{proof}
		Let us consider the optimal closed-loop dynamics, i.e., $f_{u_*} := f + b u_\ast$ with $u_\ast$ as in Theorem~\Rref{thm:optimal_feedback}. We show that $f_{u_*}$ is boundary-admissible. We begin with $u_* \in \fWpw{\infty}{w}{\Omega}$.
		By Proposition~\Rref{prop:stable_optimal_feedback}, $u_\ast$ is characterized as
		\begin{equation}
		u_\ast(z) := - \tfrac{1}{2} b(z)^\top \nabla v(z) = - b(0)^\top \hat{P} z + o(\| z \|) 
		\label{eq:lemma:exponentially_stability:feedback}
		\end{equation}
		for $z \in \Omega$ with $\| z \|$ small enough. Since $v \in \mathrm{W}^{2,\infty}(\Omega)$, we conclude $u_\ast \in \fWp{\infty}{\Omega}$ and, in particular, that $u_\ast$ is continuous. Since we have set $w = \| \cdot \|^{-1}$, \eqref{eq:lemma:exponentially_stability:feedback} together with continuity on a bounded set implies
		\[
			\| u_* \|_{\Lpw{\infty}{w}{\Omega}} \le \esssup_{z \in B_{\delta}(0)} \left|  b(0)^\top \hat{P} \tfrac{z}{\|z\|} + o(1) \right| + \delta^{-1} \esssup_{z \in \Omega \setminus B_{\delta}(0)} |u_\ast(z)| < \infty  \qquad \text{for all}\; \delta > 0
		\]
		and hence $u_* \in \fWpw{\infty}{w}{\Omega}$. Since Theorem~\Rref{thm:optimal_feedback} gives $\nu(z)^\top(f(z)+b(z)u_*(z))\le0$ for $z \in \partial \Omega$, we conclude that $u_*$ is boundary-admissible. From the above calculation we also derive
		\[
		\mathrm{D} \left( f(z) + b(z) u_\ast(z) \right)\big|_{z = 0} = \mathrm{D}f(0) + \mathrm{D} b(0) u_\ast(0) + b(0) \mathrm{D}u_\ast(0) = \mathrm{D}f(0) - b(0) b(0)^\top \hat{P},
		\]
		since $u_\ast(0) = 0$. The result from Proposition~\Rref{prop:stable_optimal_feedback} also implies that $\hat{P}$ solves the Lyapunov equation from Definition~\Rref{defi:linear_stable} \Rref{defi:stability:lqr} with respect to $f_{u_*}$ and
		\[
		\hat{H} = \hat{C}^\top\hat{C}+\hat P b(0)b(0)^\top\hat P.
		\]
		Therefore $f_{u_*}$ is even a linear stable dynamic, which
		 allows us to apply Proposition~\Rref{prop:local_uncontrolled_stability} to the closed-loop dynamic $f_{u_*}$; for any $\omega_1 > \omega_0$, there exists $\delta > 0$ such that
		\[
		\| S_\delta(t) \|_{\mathcal{L}(Y_\delta)} \le C(\omega_1) \exp( t \omega_1 ),
		\]
		where $\Omega_\delta$ and $S_\delta(t)$ are as given in Proposition~\Rref{prop:local_uncontrolled_stability} with respect to $f_{u_*}$. 
		
		Next, we show that all trajectories enter $\Omega_\delta$ within a uniform finite time. Suppose that a trajectory starting at some $z\in\Omega$ remains outside $\Omega_\delta$ throughout $[0,T]$.
		Since $\hat{P}$ is invertible~\cite[Corollary 8.4.9]{Son98},
		the norms $\| \cdot \|_{\hat{P}}$ and $\| \cdot \|$ on $\R^d$ are equivalent and it follows that there exists some $\tilde{\delta} > 0$ such that $B_{\tilde{\delta}}(0) \subseteq \Omega_\delta$ and therefore
		\[
		\Phi_{u_*}^t(z) \not \in B_{\tilde{\delta}}(0) \qquad \text{for all}\; 0 < t < T.
		\]
		By assumption, there exists an $\varepsilon > 0$ such that $\sum_{i=1}^r c_i(\Phi_{u_*}^t(z))^2 \ge \varepsilon$, and we conclude
		\[
		T \varepsilon \le \int_{0}^T \sum_{i=1}^r c_i(\Phi_{u_*}^t(z))^2 \dint t \le v(z) \le  \sup_{x \in \Omega} v(x) =: M < \infty.
		\]
		For $T := \tfrac{M + 1}{\varepsilon} < \infty$, this is a contradiction. Hence every trajectory enters $\Omega_\delta$ no later than time $T$. Since the restriction of the closed-loop dynamic to $\Omega_\delta$ is boundary-admissible, $\Omega_\delta$ is forward invariant; therefore $\Phi_{u_*}^t(z) \in \Omega_\delta$ for every $t>T$.
		We can now use this fact and the stability of $S_\delta$ to show that $S_{u_*}$ is exponentially stable. Let $\phi \in Y$. Lemma~\Rref{lemma:preadjoint_semigroup} gives
		\[
		\supp \, (S_{u_*}(T) \phi)_i \subseteq \Omega_{\delta} \qquad 0 \le i \le d.
		\]
		We denote $\fNWpw{1}{w}{\Omega_{\delta}} = Y_\delta$ and $S_\delta \colon [0,\infty) \to  \mathcal{L}(Y_\delta)$ as in Proposition~\Rref{prop:local_uncontrolled_stability}. There exists an extension
		\[
		\begin{array}{c}
			E \colon Y_\delta \to Y, \quad [(f_0, f_1, \dots, f_d)]_{\sim}  \mapsto  [ (\tilde{f}_0, \tilde{f}_1, \dots, \tilde{f}_d)]_{\sim},
			\vspace{.3cm} \\
			\text{where} \hspace{1cm}
			\tilde{f}_i(x) = \left\{\begin{array}{lcc}
				f_i(x) & \quad & \text{for}\; x \in \Omega_{\delta},\\
				0 & \quad & \text{else}.
			\end{array}\right.
		\end{array}
		\]
		The extension is well-defined on the quotient, and hence so is the map $E$.
		We conclude $\| \phi \|_{Y_\delta} = \| E \phi \|_Y$ for all $\phi \in Y_\delta$. From Proposition~\Rref{prop:local_uncontrolled_stability}, $S_\delta$ is well-defined and we have
		\[
		\| S_{u_*}(s + T) \phi \|_Y = \| S_{u_*}(s) \left( S_{u_*}(T) \phi \right) \|_Y = \| E \left( S_\delta(s) \left( S_{u_*}(T) \phi \right)\big|_{\Omega_\delta} \right) \|_Y .
		\]
		Furthermore, by the same proposition, $\| S_\delta(s) \phi\|_{Y_\delta} \le C_1 \exp(\omega_1 s) \| \phi\|_{Y_\delta}$ for all $\omega_1 > \omega_0$. By Proposition~\Rref{prop:bound_semigroup}, we have $\|S_{u_*}(T)\|_{\mathcal{L}(Y)} \le C_2 < \infty$, hence
		\[
		\| S_{u_*}(s + T) \phi \|_Y \le C_1 C_2 \exp(\omega_1 s )  \| \phi \|_Y \quad \text{for all}\; 0 \le s < \infty
		\]
		or equivalently
		\[
		\| S_{u_*}(t) \phi \|_Y \le C \exp(\omega_1 t )  \| \phi \|_Y \qquad \text{for} \; 0 \le t < \infty,
		\]
		where we define 
		\[
		C = \max \left\{ C_1 C_2 ,\underset{ 0 \le t \le T}{\sup} \, \left\{ \exp(-\omega_1 t) \| S_{u_*}(t) \|_{\mathcal{L}(Y)} \right\} \right \}.
		\]
		It remains to verify the assertion on $X$. Set $X_\delta:=\fLpw{1}{w}{\Omega_\delta}$. The first estimate in the proof of Proposition~\Rref{prop:local_uncontrolled_stability} gives
		\[
		\|S_\delta(s)q\|_{X_\delta}\leq C_3\exp(\omega_1s)\|q\|_{X_\delta}.
		\]
		Zero extension $X_\delta\to X$ is isometric and restriction $X\to X_\delta$ is contractive. The same support argument at time $T$ therefore yields
		\[
		\|S_{u_*}(s+T)\phi\|_X
		\leq C_3\exp(\omega_1s)\|S_{u_*}(T)\phi\|_X
		\leq C_3C_4\exp(\omega_1s)\|\phi\|_X,
		\]
		where $C_4=\|S_{u_*}(T)\|_{\mathcal L(X)}<\infty$. Enlarging the constant on $[0,T]$ proves exponential stability on $X$ with the same admissible types $\omega_1>\omega_0$.
	\end{proof}
	\begin{rmrk}
		Assumptions \Rref{req:exp_stab:rank} and \Rref{req:exp_stab:positivity} of the lemma can be weakened by a slight modification of the weighted Sobolev space. In particular, a weight other than $w(x) = \sfrac{1}{\| x \|}$ must be chosen; this is not pursued in this manuscript.
		Note that, for the case of a full linear observation, $c_i(x) := x_i$, both assumptions are trivially satisfied.
	\end{rmrk}
	Next, we introduce the value bilinear form for the control problem associated with \eqref{eq:dynamics_intro}.
	\begin{dfntn}
		\label{defi:value_bilinear_form}
		Let $(f,b,u)$ be a boundary-admissible closed-loop dynamic such that $S_u$ is exponentially stable over $Y$, and let $\{c_i\}_{i \in \N} \subseteq Y^\ast$ be a nuclear sequence of observables with $C_u$ the corresponding observation with control penalty. We define the \emph{value bilinear form} as
		\[
		\begin{array}{rlcl}
			\einner{P_u} \colon & Y \times Y & \to & \R, \\
			& (\phi,\psi) & \mapsto & \int_0^\infty \inner[\ell_2]{C_{u} S_{u}(t) \phi}{C_{u} S_{u}(t)\psi}\dint{t}.
		\end{array}
		\]
	\end{dfntn}
	\begin{lmm}
		\label{lem:lowrank_control}
		Let $(f,b,u)$ be a boundary-admissible closed-loop dynamic such that $S_u$ is exponentially stable over $Y$, and let $\{c_i\}_{i \in \N} \subseteq Y^\ast$ be a nuclear sequence of observables with $C_u$ the corresponding observation with control penalty. Then:
		\begin{enumerate}[label=(\roman*)]
			\item The corresponding value bilinear form $\einner{P_{u}}$ admits the decomposition
			\begin{equation}
				\inner[P_{u}]{\phi}{\psi} = \sum_{i=1}^\infty \sigma_i \inner[Y,Y^\ast]{\phi}{p_i} \inner[Y,Y^\ast]{\psi}{p_i} \qquad \text{for all}\; \phi,\psi \in Y
				\label{lemma:lowrank_decomposition:eq:decomposition}
			\end{equation}
			for some $p_i \in Y^\ast$ with $\|p_i\|_{H_Y} =1$, $\{\sigma_i\}_{i \in \N} \subseteq \R_+ \setminus \{0\}$, and $\sum_{i=1}^\infty \sigma_i < \infty$.
			\label{lem:lowrank_control:decomposition} 		
		\end{enumerate}
	\end{lmm}
	\begin{proof}
		In \cite[Theorem 3.4]{BreH23} a decomposition of the form \eqref{lemma:lowrank_decomposition:eq:decomposition} is obtained by fixing an orthonormal basis $(H_n)$ in time and defining $p_{ni}\in H_X$ via the Riesz representation theorem by
		\[
		\inner[H_X]{\phi}{p_{ni}}
		= \int_0^\infty \inner[H_X]{S(t)\phi}{c_i}\, H_n(t)\, \mathrm{d}t .
		\]
		That proof assumes exponential stability of $S$ on $X$. In our setting, $S_u$ is exponentially stable on $Y$, and by Lemma~\Rref{lemma:gelfand_like} the pair $(Y,H_Y)$ satisfies the same properties required for $X$ and $H_X$ in \cite{BreH23}. Replacing $X,H_X$ with $Y,H_Y$ then yields \Rref{lem:lowrank_control:decomposition}.
		
	\end{proof}
	\begin{thrm}
		\label{thm:lowrank_optimalcontrol}
		Let $(f,b)$ be a boundary-admissible controlled dynamic pair, $\{c_i\}_{i \in \N}$ a nuclear sequence of observables, and let $u_*$ denote the corresponding optimal feedback from Theorem~\Rref{thm:optimal_feedback}. If $u_* \in Y^\ast$ is boundary-admissible and $S_{u_*}$ is exponentially stable over $Y$, then:
		\begin{enumerate}[label=(\roman*)]
			\item The value function from Definition~\Rref{defi:value_function} satisfies $v\in \fWpw{1}{w^2}{\Omega}$ and $v(z) = \sum_{i=1}^\infty \sigma_i p_i(z)^2$ for all $z \in \Omega$.
			\label{thm:lowrank_optimalcontrol:value_function}
			\item \label{thm:lowrank_optimalcontrol:feedback} $u_*(z) = - \sum_{i=1}^\infty \sigma_i (b(z)^\top \nabla p_i(z)) p_i(z)$ for a.e. $z \in \Omega$, with convergence in $\fLpw{1}{w}{\Omega}$.	 		
		\end{enumerate}
		Here $\{\sigma_i\}_{i\in\N}$ and $\{ p_i\}_{i\in\N}$ denote the decomposition from Lemma~\Rref{lem:lowrank_control}.
	\end{thrm}
	\begin{proof}
		The proof of a similar sum-of-squares representation from \cite[Theorem 3.10]{BreH23} never uses orthogonality with respect to $H_X$, so the result holds for any decomposition of the form \eqref{lemma:lowrank_decomposition:eq:decomposition}, which shows \Rref{thm:lowrank_optimalcontrol:value_function}. Finally, \Rref{thm:lowrank_optimalcontrol:feedback} follows directly from Theorem~\Rref{thm:optimal_feedback}.
	\end{proof}
	\begin{rmrk}
		Note that the previous results hold, in particular, for systems that satisfy the assumptions of Lemma~\Rref{lemma:exponentially_stability}.
	\end{rmrk}
	
	\begin{lmm}
		\label{lemma:decay}
		Let $(f,b,u)$ be a boundary-admissible closed-loop dynamic such that $S_u$ is exponentially stable over $Y$, and let $\{c_i\}_{i \in \N} \subseteq Y^\ast$ be a nuclear sequence of observables with $C_u$ the corresponding observation with control penalty. Let $m \in \N$ be even and suppose that $c_i\equiv 0$ for all $i>r$ and some $r\in\N$, that $(f+bu)_j\in W^{m,\infty}(\Omega)$ for $1\leq j\leq d$, and that $u,c_i\in W^{m+1,\infty}(\Omega)$ for $1\leq i\leq r$.
		
		Then the decomposition from Lemma~\Rref{lem:lowrank_control} can be chosen such that
		\[
		\sum_{i=N}^\infty \sigma_i = O(N^{-m}) \qquad \text{as }N\to\infty.
		\]
	\end{lmm}
	\begin{proof}
		The stated regularity implies that $u,c_1,\dots,c_r$ belong to the domains of the first $m$ powers of the weak-$\ast$ generator of $S_u^\ast$ on $Y^\ast$. The assertion thus follows from the proof of \cite[Lemma 3.7]{BreH23}, applied on $Y$.
	\end{proof}
	The following auxiliary lemma shows that the eigenfunctions of the value bilinear form span an invariant subspace for $S_u^\ast$ in $H_Y$.  
	\begin{lmm}
		\label{lemma:invariance}
		Let $(f,b,u)$ be a boundary-admissible closed-loop dynamic such that the semigroup $S_u(t)$ is exponentially stable over $Y$, and let $\{c_i\}_{i \in \N}$ be a nuclear sequence of observables. Furthermore, let $p_i,\sigma_i$ denote a decomposition as in Lemma~\Rref{lem:lowrank_control}. Then the following hold:
		\begin{enumerate}[label=(\roman*)]
			\item
			\label{lemma:invariance:invariance}
			The space $\clSpan \{ p_1, p_2, \dots\} \subseteq H_Y$
			is invariant under $S_u^\ast(t)$ for all $0 \le t < \infty$.
			\item $c_i \in \clSpan \{ p_1, p_2, \dots \} \subseteq H_Y$ for all $i\in\N_0$, where $c_0:=u$.
			\label{lemma:invariance:element}
		\end{enumerate}
	\end{lmm}
	\begin{proof}
		First, we show that $\clSpan \{ p_1, p_2, \dots \}^\perp$ is $S_u(t)$-invariant. \\
		Let $\phi \in \clSpan \{ p_1, p_2, \dots \}^{\perp}$. Since $\inner[H_Y]{\phi}{p_i} = 0$, it follows that 
		\[
			\inner[P_{u}]{\phi}{\phi} = \sum_{j=1}^\infty \sigma_j \inner[H_Y]{\phi}{p_j}^2 = 0,
		\] 
		which leads to
		\begin{align*}
			0 \le & \sum_{i=1}^\infty \sigma_{i} \inner[H_Y]{S_u(t)\phi}{p_i}^2 = \inner[P_u]{S_u(t)\phi}{S_u(t)\phi}\\
			= &\int_0^\infty \inner[\ell_2]{C_u S_u(s) S_u(t) \phi}{C_u S_u(s) S_u(t) \phi} \dint s \\
			= &\int_0^\infty \inner[\ell_2]{C_u S_u(t+s) \phi}{C_u S_u(t+s) \phi} \dint s \\
			= & \int_0^\infty \inner[\ell_2]{C_u S_u(s)\phi}{C_u S_u(s) \phi} \dint s - \int_{0}^t \inner[\ell_2]{C_u S_u(s)\phi}{C_u S_u(s) \phi} \dint s\\
			= & \inner[P_u]{\phi}{\phi} - \int_{0}^t \inner[\ell_2]{C_u S_u(s)\phi}{C_u S_u(s) \phi} \dint s \le 0.
		\end{align*}
		Therefore $\inner[H_Y]{S_u(t)\phi}{p_i} = 0$. We conclude
		\[
		\inner[H_Y]{\phi}{S_u^\ast(t) p_i}  = \inner[H_Y]{S_u(t) \phi}{p_i} = 0 \qquad \text{for all}\; \phi \in \clSpan \{ p_1, p_2, \dots \}^{\perp}
		\]
		and it follows that $S_{u}^\ast(t) p_i \in \left( \clSpan \{ p_1, p_2, \dots \}^\perp \right)^\perp \subseteq \clSpan \{ p_1, p_2, \dots \}$.\\
		Now let us show \Rref{lemma:invariance:element}:
		Let $\phi \in \clSpan \{ p_1, p_2, \dots \}^{\perp}$. We have
		\begin{align*}
			0 \le & \int_0^t \inner[\ell_2]{C_u S_u(s)\phi}{C_u S_u(s) \phi} \dint s \\
			= & \int_0^\infty \inner[\ell_2]{C_u S_u(s)\phi}{C_u S_u(s) \phi} \dint s - \int_{t}^\infty \inner[\ell_2]{C_u S_u(s)\phi}{C_u S_u(s) \phi} \dint s\\
			= & \inner[P_u]{\phi}{\phi} - \int_{t}^\infty \inner[\ell_2]{C_u S_u(s)\phi}{C_u S_u(s) \phi} \dint s \le 0.
		\end{align*}
		By the fundamental theorem of calculus \cite[Theorem E3.6 (2)]{Alt12} we conclude
		\[
		\inner[\ell_2]{C_u \phi}{C_u \phi} = \lim_{t \searrow 0} \frac{1}{t} \int_0^t \inner[\ell_2]{C_u S_u(s)\phi}{C_u S_u(s) \phi} \dint s = 0
		\]
		and therefore
		$
		\sum_{i=0}^\infty \inner[H_Y]{c_i}{\phi}^2  = 0.
		$
		We conclude 
		\[
		c_i \in \left( \clSpan\{p_1, p_2, \dots \}^\perp\right)^\perp \subseteq \clSpan\{p_1, p_2, \dots \}.
		\]
	\end{proof}
	Let us formulate the Lyapunov equation from \cite{BreH23} in a weaker form.
	\begin{lmm}
		\label{lemma:lyapunov}
		Let $(f,b,u)$ be a boundary-admissible closed-loop dynamic such that
		$S_u(t)$ is exponentially stable on $Y$, and let $\{c_i\}_{i \in \N}$ be a nuclear sequence of observables.
		Let $\einner{P_u}$ be as in Definition~\Rref{defi:value_bilinear_form}, and let $\{p_k\}_{k\in\N}$ and $\{\sigma_k\}_{k\in\N}$ be the decomposition from Lemma~\Rref{lem:lowrank_control}. Then
		\[
			0 = \sum_{k=1}^\infty \sigma_k \left( \inner[H_X]{\phi}{ A_u^\ast p_k} \inner[H_X]{\psi}{p_k} + \inner[H_X]{\phi}{p_k} \inner[H_X]{\psi}{A_u^\ast p_k} \right) + \inner[\ell_2]{C_u  \phi}{C_u \psi} \qquad \text{for all}\; \phi,\psi \in H_X.
		\]
		Here $A^\ast_u \in \mathcal{L}(H_Y;H_X)$ is defined by
		\[
			A_u^\ast \colon H_Y \to H_X, \quad \phi \to (f + b u)^\top \nabla \phi.
		\]
	\end{lmm}
	\begin{proof}
		Since every component of $f_u$ belongs to $X^\ast$, for every $p\in H_Y$ we have
		\[
		\|A_u^\ast p\|_{H_X}
		=\|f_u^\top\nabla p\|_{H_X}
		\leq d\max_{1\leq j\leq d}\|(f_u)_j\|_{X^\ast}\|p\|_{H_Y}.
		\]
		Thus $A_u^\ast\in\mathcal L(H_Y;H_X)$. Moreover, Proposition~\Rref{prop:generator_koopman}, applied to the closed-loop dynamic $f_u := f +b u$, shows that $H_Y$ is contained in the domain of the generator of the Koopman semigroup on $H_X$ and that the restriction of this generator to $H_Y$ is $A_u^\ast$.
		
		Fix $\phi,\psi\in H_X$ and identify them with $[(\phi,0,\dots,0)]_\sim,[(\psi,0,\dots,0)]_\sim\in H_Y$. The explicit preadjoint formula from Definition~\Rref{defi:preadjoint_semigroup} shows that the Perron--Frobenius semigroups on $H_X$ and $H_Y$ agree under this identification. Moreover, the dual pairing with $p_k\in H_Y\subset H_X$ agrees with the $H_X$ inner product. Hence Lemma~\Rref{lem:lowrank_control} gives
		\[
		\inner[P_u]{\phi}{\psi}
		=\sum_{k=1}^\infty\sigma_k
		\inner[H_X]{\phi}{p_k}\inner[H_X]{\psi}{p_k}.
		\]
		The series is absolutely convergent because $\|p_k\|_{H_Y}=1$ and $\sum_{k=1}^\infty\sigma_k<\infty$. For $h\geq0$, the semigroup property and Definition~\Rref{defi:value_bilinear_form} yield
		\begin{align*}
			\inner[P_u]{S_u(h)\phi}{S_u(h)\psi}
			=\int_h^\infty
			\inner[\ell_2]{C_uS_u(t)\phi}{C_uS_u(t)\psi}\dint t
			=\inner[P_u]{\phi}{\psi}
			-\int_0^h
			\inner[\ell_2]{C_uS_u(t)\phi}{C_uS_u(t)\psi}\dint t.
		\end{align*}
		Strong continuity of $S_u$ and boundedness of $C_u$ therefore imply
		\begin{equation}
		\lim_{h\searrow0}\frac{1}{h}
		\left(\inner[P_u]{S_u(h)\phi}{S_u(h)\psi}
		-\inner[P_u]{\phi}{\psi}\right)
		=-\inner[\ell_2]{C_u\phi}{C_u\psi}.
		\label{eq:lemma:lyapunov:CCT_term}
		\end{equation}
		On the other hand, for every $k\in\N$, the generator identity above gives
		\begin{equation}
		\lim_{h\searrow0}\frac{1}{h} \left(\inner[H_X]{S_u(h)\phi}{p_k} -\inner[H_X]{\phi}{p_k}\right)
		= \lim_{h\searrow0}\frac{1}{h} \inner[H_X]{\phi}{S^\ast_u(h)p_k - p_k} 
		= \inner[H_X]{\phi}{A_u^\ast p_k}.
		\label{eq:lemma:lyapunov:diff_quotient}
		\end{equation}
		For a fixed $h_0>0$, Proposition~\Rref{prop:bound_semigroup} gives
		\[
		M:=\sup_{0\leq h\leq h_0}
		\|S_u^\ast(h)\|_{\mathcal L(H_X)}<\infty
		\quad\text{and}\quad
		K:=\|A_u^\ast\|_{\mathcal L(H_Y;H_X)}.
		\]
		Since $p_k \in H_Y\subseteq\mathcal D(A_u^\ast)$, the standard generator identity gives
		\[
			S_u^\ast(h)p_k-p_k=\int_0^h S_u^\ast(s)A_u^\ast p_k\,\dint s.
		\]
		Therefore, for $0<h\leq h_0$,
		\[
		\left\|\tfrac{1}{h} (S_u^\ast(h)p_k-p_k)\right\|_{H_X}
		\leq M\|A_u^\ast p_k\|_{H_X}\leq MK
		\quad \text{and} \quad 
		\|S_u^\ast(h)p_k\|_{H_X}\leq M.
		\]
		Consequently, the corresponding difference quotient of the product is bounded by $KM(M+1)\|\phi\|_{H_X}\|\psi\|_{H_X}$ independently of $k$ and $h$. Since $\sum_k\sigma_k<\infty$, dominated convergence for series yields
		\begin{align*}
			&\lim_{h\searrow0}\frac{1}{h}
			\left(\inner[P_u]{S_u(h)\phi}{S_u(h)\psi}
			-\inner[P_u]{\phi}{\psi}\right)\\
			&\qquad=\sum_{k=1}^\infty\sigma_k\left(
			\inner[H_X]{\phi}{A_u^\ast p_k}\inner[H_X]{\psi}{p_k}
			+\inner[H_X]{\phi}{p_k}\inner[H_X]{\psi}{A_u^\ast p_k}
			\right).
		\end{align*}
		Together with \eqref{eq:lemma:lyapunov:CCT_term} this identity proves the assertion.
	\end{proof}
	We now come to one of our main contributions, which shows that, through a coordinate transform, the value function and the optimal feedback can be linked to a quadratic operator equation on $\ell_2$.
	To prepare for the theorem in which the operator equation is stated, we first give auxiliary propositions and a definition introducing all quantities that will appear in it.
	\begin{dfntn}
		\label{defi:auxiliary_nonlinear_op_eq}
		Assume that $(f,b)$ is a boundary-admissible controlled dynamic pair, that
		$\{c_i\}_{i\in\N}$ is a nuclear sequence of observables, and that the
		Perron--Frobenius semigroup $S_{u_*}(t)$ is exponentially stable on $Y$.
		Let $u_*$ be the optimal feedback from Theorem~\Rref{thm:optimal_feedback},
		and let $\{p_i\}_{i\in\N}$ and $\{\sigma_i\}_{i\in\N}$ be the corresponding
		decomposition from Lemma~\Rref{lem:lowrank_control}.
		
		Let $\{v_i\}_{i\in\N}\subseteq Y^\ast$ be a Riesz sequence in $H_Y$ such that
		\[
		\clSpan\{p_i:i\in\N\}\subseteq
		V:=\clSpan\{v_i:i\in\N\}\subseteq H_Y,
		\]
		and denote its synthesis operator by
		\[
		T:\ell_2\to V,
		\qquad
		Ta=\sum_{i=1}^\infty a_i v_i.
		\]
		Let $\mathbb P_V:H_Y\to V$ be the $H_Y$-orthogonal projection. Moreover,
		if $\Psi^{-1}:\fNWpw{2}{w}{\Omega}\to\fWpw{2}{w}{\Omega}$ is the isometry
		from Lemma~\Rref{lemma:iso_hilbert}, define
		\[
		E:\fLpw{2}{w}{\Omega}\to\fWpw{2}{w}{\Omega},
		\qquad
		E\phi:=\Psi^{-1}\bigl([(\phi,0,\ldots)]_{\sim}\bigr).
		\]
		
		Define the bounded linear operators
		$M,F^*,H,F_k^*\in\mathcal L(\ell_2)$, $k\in\N$, by
		\[
		\begin{aligned}
			M a &:= \bigl(\inner[H_X]{v_i}{Ta}\bigr)_{i\in\N},&
			F^*a &:= \bigl(\inner[H_X]{v_i}{f^\top\nabla(Ta)}\bigr)_{i\in\N},\\
			F_k^*a &:= \bigl(\inner[H_X]{v_i}{v_k b^\top\nabla(Ta)}\bigr)_{i\in\N},&
			H a &:= \bigl(\inner[H_X]{c_i}{Ta}\bigr)_{i\in\N},
		\end{aligned}
		\]
		and set $F:=(F^*)^*$ and $F_k:=(F_k^*)^*$. For every $a\in\ell_2$ such that $Ta\in X^\ast$, also define
		\[
			F[a]^*c:=\bigl(\inner[H_X]{v_i}{(Ta)b^\top\nabla(Tc)}\bigr)_{i\in\N},
			\qquad c\in\ell_2,
		\]
		and $F[a]:=(F[a]^*)^*$. In particular, $F[e_k]=F_k$ for the canonical unit vectors $e_k\in\ell_2$. For
		\[
		\widetilde\Sigma =\sum_{k=1}^\infty\widetilde\sigma_k\inner[\ell_2]{\cdot}{a_k}a_k\in\mathcal N_{\operatorname{sa}}(\ell_2),
		\]
		where the displayed series is the spectral decomposition of $\widetilde\Sigma$ and the vectors $a_k$ are orthonormal, define
		\[
		\beta(\widetilde\Sigma):=-\left(\sum_{k=1}^\infty\widetilde\sigma_k\inner[H_X]{v_i b^\top\nabla(Ta_k)}{Ta_k} \right)_{i\in\N},
		\]
		and
		\[
		\mathcal{D}(\alpha) := \left\{ \widetilde\Sigma\in\mathcal N_{\operatorname{sa}}(\ell_2): \beta(\widetilde\Sigma)\in\Image M,\quad
		T M^{-1}\beta(\widetilde\Sigma)\in X^\ast \right\}.
		\]
		Since $M$ is injective by Proposition~\Rref{prop:nonlinear_op_eq_bd}, the map
		\[
		\alpha:\mathcal D(\alpha)\to\ell_2, \qquad \alpha(\widetilde\Sigma):=M^{-1}\beta(\widetilde\Sigma)
		\]
		is well defined. Finally, for $\Sigma\in\mathcal D(\alpha)$ and $a,b\in\Image M$, define
		\[
		\begin{aligned}
			\inner[G(\Sigma)]{a}{b}:={}&
			\inner[\Sigma]{F[\alpha(\Sigma)]M^{-1}a}{b}+\inner[\Sigma]{a}{F[\alpha(\Sigma)]M^{-1}b}\\
			& \qquad +\inner[\ell_2]{a}{\alpha(\Sigma)} \inner[\ell_2]{b}{\alpha(\Sigma)}.
		\end{aligned}
		\]
	\end{dfntn}
	
	The following proposition shows that the quantities in the definition are well-defined. 
	\begin{prpstn}
		\label{prop:nonlinear_op_eq_bd}
		Under the assumptions and definitions of Definition~\Rref{defi:auxiliary_nonlinear_op_eq}, the following hold:
		\begin{enumerate}[label=(\roman*)]
			\item $E$ is bounded and injective.
			\label{item:prop:nonlinear_op_eq_bd:E}
			\item $M = T^\ast E T$ is bounded, injective and self-adjoint.
			\label{item:prop:nonlinear_op_eq_bd:M}
			\item For any $\tilde{f} \colon \Omega \to \R^d$ with $\tilde{f}_i \in X^\ast$ for $1 \le i \le d$, the operator
			\[
			\tilde{F}^\ast \colon \ell_2 \to \ell_2, \qquad a \mapsto \left( \inner[H_X]{v_i}{\tilde{f}^\top \nabla ( T a)} \right)_{i \in \N}
			\]
			is bounded, i.e., $\tilde{F}^\ast \in \mathcal{L}(\ell_2)$. In particular, $F^\ast$, $F^\ast_k$, and $F[a]^\ast$ for every $a\in\ell_2$ with $Ta\in X^\ast$ are bounded.
			\label{item:prop:nonlinear_op_eq_bd:F}
			\item $H$ is bounded.
			\label{item:prop:nonlinear_op_eq_bd:H}
		\end{enumerate}
	\end{prpstn}
	\begin{proof}
		See Appendix~\Rref{sec:technical_proofs}.
	\end{proof}
	\begin{thrm}
		\label{thm:nonlinear_op_eq}
		Under the assumptions and definitions of Definition~\Rref{defi:auxiliary_nonlinear_op_eq} there exists a positive, nuclear linear operator $\Sigma$ with eigenvectors $\{a_k\}_{k \in \N}$ and eigenvalues $\{\sigma_k\}_{k \in \N} \subseteq \R_+\setminus \{0\}$, i.e.,
		\[
		\Sigma := \sum_{k=1}^\infty \sigma_k \inner[\ell_2]{\cdot}{a_k} a_k \in \mathcal{D}(\alpha) \subseteq \mathcal{N}_{\operatorname{sa}}(\ell_2),
		\] 
		such that the following statements hold:
		\begin{enumerate}[label=(\roman*)]
			\item The value bilinear form satisfies 
			\[
			\inner[P_{u_*}]{\phi}{\psi} = \inner[\Sigma]{ T^\ast \mathbb{P}_V \phi }{ T^\ast \mathbb{P}_V \psi }  \qquad \text{for all} \;  \phi,\psi \in H_Y.
			\]
			\item For almost every $z \in \Omega$,
			\[
			v(z) = \sum_{k=1}^\infty \sigma_k ((T a_k)(z))^2  \quad \text{and} \quad u_*(z) = - \sum_{k=1}^\infty \sigma_k (T a_k)(z) b^\top \nabla (T a_k)(z).
			\]
			\item For all $a,b \in \Image M$,
			\begin{equation*}
				\inner[\Sigma]{F M^{-1} a}{b} + \inner[\Sigma]{a}{F M^{-1} b} + \inner[G(\Sigma)]{a}{b} + \inner[\ell_2]{H M^{-1}a}{H M^{-1}b} = 0.
			\end{equation*}
			\label{eq:nonlinear_op_eq:non_linear_operator_eq}
		\end{enumerate}
	\end{thrm}
	\begin{proof}
		The idea of the proof is to apply a coordinate transform via $T$ to the Lyapunov equation from Lemma~\Rref{lemma:lyapunov}. We start by defining $\Sigma$ via the coordinate transform of $\einner{P_{u_*}}$ under $T$, i.e.,
		\[
			\inner[\Sigma]{a}{b} := \inner[P_{u_*}]{T^{-\ast} a}{T^{-\ast} b} \qquad \text{for all}\; a,b \in \ell_2.
		\]
		Since $T^{-\ast}$ is bounded and $\einner{P_{u_*}}$ admits the low-rank decomposition from Lemma~\Rref{lem:lowrank_control} with summable eigenvalues, $\Sigma$ is positive and nuclear. Let $\phi,\psi \in H_Y$. Since $p_i\in V$ for every $i\in\N$, we have
		\begin{align*}
			\inner[P_{u_*}]{\phi}{\psi}
			=\sum_{i=1}^\infty \sigma_i \inner[H_Y]{\phi}{p_i}\inner[H_Y]{\psi}{p_i}
			=\sum_{i=1}^\infty \sigma_i \inner[H_Y]{\mathbb P_V\phi}{p_i}\inner[H_Y]{\mathbb P_V\psi}{p_i}
			=\inner[\Sigma]{T^\ast\mathbb P_V\phi}{T^\ast\mathbb P_V\psi}.
		\end{align*}
		This proves (i). Next, we identify $A_{u_*}$ and $C_{u_*}$ in the new coordinates. By Lemma~\Rref{lemma:invariance}, $u_* \in V$ and therefore $T^{-1}u_*\in\ell_2$. Since $\Sigma$ is positive and nuclear, it has eigenvalues $\{\sigma_k\}_{k\in\N}$ and orthonormal eigenvectors $\{a_k\}_{k\in\N}$ such that
		\begin{equation}
			\inner[P_{u_*}]{\phi}{\psi} = \inner[\Sigma]{T^{\ast} \mathbb P_V\phi}{T^{\ast} \mathbb P_V\psi} = \sum_{k=1}^\infty \sigma_k \inner[H_Y]{\phi}{T a_k}\inner[H_Y]{\psi}{T a_k} \qquad \text{for all}\; \phi,\psi \in H_Y.
			\label{eq:lowrank_Sigma_decomposition}
		\end{equation}
		As stated in the proof of Theorem~\Rref{thm:lowrank_optimalcontrol}, the argument in \cite[Theorem 3.10]{BreH23} only requires a decomposition with summable coefficients. Therefore, Theorem~\Rref{thm:optimal_feedback} yields
		\[
			v(z) = \sum_{k=1}^\infty \sigma_k (T a_k)(z)^2 \quad \text{and} \quad u_*(z) = -\sum_{k=1}^\infty \sigma_k b(z)^\top \nabla (T a_k)(z)(T a_k)(z)
		\]
		for almost every $z \in \Omega$. Using $M=T^\ast ET$ and the convergence of the feedback series in $X$, we obtain
		\[
			M T^{-1} u_* = T^{\ast} E u_* = \left( \inner[H_X]{v_k}{u_*} \right)_{k \in \N} = \beta(\Sigma).
		\]
		Thus $\beta(\Sigma) \in \Image M$ and $M^{-1}\beta(\Sigma)=T^{-1}u_*$. Since $u_*\in Y^\ast\subseteq X^\ast$, we have $T M^{-1}\beta(\Sigma)=u_*\in X^\ast$. Hence $\Sigma\in\mathcal D(\alpha)$ and $\alpha(\Sigma)=T^{-1}u_*$, which yields
		\begin{equation}
			u_* = T \alpha(\Sigma) = \sum_{k=1}^\infty \alpha(\Sigma)_k v_k.
			\label{eq:nonlinear_op_eq:u_repr}
		\end{equation}
		Proposition~\Rref{prop:nonlinear_op_eq_bd} shows that
		$
			F + F[\alpha(\Sigma)] \in \mathcal L(\ell_2)
		$. Moreover,
		\[
			\sum_{k=1}^\infty \sigma_k \|T a_k\|_{H_Y}^2 \leq \|T\|^2 \sum_{k=1}^\infty \sigma_k < \infty.
		\]
		Consequently, the proof of Lemma~\Rref{lemma:lyapunov} applies to the decomposition in \eqref{eq:lowrank_Sigma_decomposition} and gives, for all $\phi,\psi\in H_X$,
		\begin{equation}
			0 = \sum_{k=1}^\infty \sigma_k \left( \inner[H_X]{\phi}{A_{u_*}^\ast T a_k} \inner[H_X]{\psi}{T a_k} + \inner[H_X]{\phi}{T a_k} \inner[H_X]{\psi}{A_{u_*}^\ast T a_k} \right) + \inner[\ell_2]{C_{u_*}\phi}{C_{u_*}\psi}.
			\label{eq:thm:operator_equation:weak_lyapunov}
		\end{equation}
		For $\phi = T M^{-1} a$ and $\psi = T M^{-1} b$, where $a,b\in\Image M$, we have
		\[
			\inner[H_X]{T M^{-1} a}{T a_k} = \inner[H_Y]{E T M^{-1} a}{T a_k} = \inner[\ell_2]{a}{a_k}.
		\]
		By the definitions of $F^\ast$ and $F[\alpha(\Sigma)]^\ast$, we also have
		\[
			\inner[H_X]{T M^{-1} a}{A_{u_*}^\ast T a_k} = \inner[\ell_2]{(F + F[\alpha(\Sigma)]) M^{-1} a}{a_k}.
		\]
		In addition, Definition~\Rref{defi:control_operator} and \eqref{eq:nonlinear_op_eq:u_repr} yield
		\begin{align*}
			\inner[\ell_2]{C_{u_*} T M^{-1} a}{C_{u_*} T M^{-1} b}
			=\inner[\ell_2]{H M^{-1}a}{H M^{-1}b} +\inner[\ell_2]{a}{\alpha(\Sigma)}\inner[\ell_2]{b}{\alpha(\Sigma)}.
		\end{align*}
		Testing \eqref{eq:thm:operator_equation:weak_lyapunov} with $\phi = T M^{-1} a$ and $\psi = T M^{-1}b$ gives
		\begin{align*}
			0={}&\sum_{k=1}^\infty \sigma_k \left( \inner[\ell_2]{(F + F[\alpha(\Sigma)]) M^{-1} a}{a_k} \inner[\ell_2]{b}{a_k} + \inner[\ell_2]{a}{a_k} \inner[\ell_2]{(F + F[\alpha(\Sigma)]) M^{-1}b}{a_k} \right)\\
			& \qquad +\inner[\ell_2]{H M^{-1} a}{H M^{-1} b} + \inner[\ell_2]{a}{\alpha(\Sigma)}\inner[\ell_2]{b}{\alpha(\Sigma)}.
		\end{align*}
		Rearranging the terms involving $F[\alpha(\Sigma)]$ and using the definition of $G(\Sigma)$ proves~\Rref{eq:nonlinear_op_eq:non_linear_operator_eq}.
	\end{proof}
	\begin{rmrk}
		\label{rem:riccati}
		In the LQR case, the operator equation 
		from Theorem~\Rref{thm:nonlinear_op_eq}\ref{eq:nonlinear_op_eq:non_linear_operator_eq} 
		reduces exactly to the algebraic Riccati equation. To see this, consider
		\[
		f(x) := \hat{A} x, \quad b(x) := \hat{b}, \quad \text{and} \quad c_i(x) := \hat{c}_i^\top x, \quad 1\le i\le r,
		\]
		and set $\hat C:=(\hat c_1,\ldots,\hat c_r)^\top$.
		Then, for any boundary-admissible linear feedback $u$, the eigenfunctions $p_i$ of $\einner{P_u}$ are linear functions. In the corresponding finite-dimensional setting, we can therefore choose $v_i(x) := x_i$, $1\le i\le d$, as a Riesz basis for the space of linear observables. Furthermore, define
		\[
		M_{ij} = \int_\Omega  x_i x_j \;w^2(x) \dint x \quad \text{and} \quad \Gamma_{ijk} := \inner[H_X]{v_i b^\top \nabla v_j}{v_k} = \hat{b}_j M_{ik}.
		\]
		The quantities in Definition~\Rref{defi:auxiliary_nonlinear_op_eq} can then be identified as
		\begin{align*}
			F^\ast_{ij} = 
			\int_\Omega x_i x^\top \hat{A}^\top e_j\;  w^2(x) \dint x 
			= (M \hat{A}^\top)_{ij},
			\quad 
			(F_k^\ast)_{ij} = \Gamma_{ijk}, \quad H = \hat C M, \\
			\beta(\Sigma)_k = -\sum_{i,j=1}^d \Sigma_{ij} \Gamma_{ijk} = -\sum_{i,j = 1}^d \Sigma_{ij} \hat{b}_j M_{ik}, \quad \text{and} \quad \alpha(\Sigma) = -\Sigma\hat{b}.
		\end{align*}
		A simple calculation shows that $G(\Sigma)=-\Sigma\hat{b}\hat{b}^\top\Sigma$, so the operator equation reduces to the algebraic Riccati equation
		\[
			0=\hat{A}^\top\Sigma+\Sigma\hat{A}-\Sigma\hat{b}\hat{b}^\top\Sigma+\hat C^\top\hat C.
		\]
		This equation has a unique symmetric positive semidefinite stabilizing solution $\Sigma_{\mathrm{LQR}}$ under standard assumptions (see, e.g., \cite{BenS13}). For the corresponding unconstrained LQR problem on $\mathbb R^d$, the value function and the unique optimal feedback are
		\[
			v_{\mathbb R^d}(z)=z^\top\Sigma_{\mathrm{LQR}}z
			\qquad\text{and}\qquad
			u_{\mathrm{LQR}}(z)=-\hat b^\top\Sigma_{\mathrm{LQR}}z,
		\]
		respectively. Therefore, if we assume that $u_{\mathrm{LQR}}$ is boundary-admissible on $\Omega$, then $v=v_{\mathbb R^d}$ and $u_*=u_{\mathrm{LQR}}$ on $\Omega$. In this case, the state-constrained optimal feedback is linear, and the operator equation from Theorem~\Rref{thm:nonlinear_op_eq}\ref{eq:nonlinear_op_eq:non_linear_operator_eq} is precisely the algebraic Riccati equation.
	\end{rmrk}
	\section{Numerical proof of concept}
	\label{sec:numerics}
	
	In this section, we present a short numerical example in two dimensions to illustrate the theoretical findings from Section~\Rref{ssec:koopman_feedback}. In particular, we investigate the decay of the eigenvalues when working in the space~$H_Y$.
	
	We use a discretization technique similar to that in \cite{BreH23} and choose tensor products of splines of degree $p=5$ with $n=31$ grid points, orthonormalized with respect to~$H_Y$, as the Riesz sequence $\{v_i\}_{i\in\N}$ in Theorem~\Rref{thm:nonlinear_op_eq}. To approximate the inner products with respect to~$H_X$, we employ a quadrature rule based on Gauss--Legendre integration on each subinterval. Since we use only finitely many basis elements, the operator equation over~$\ell_2$ becomes a quadratic matrix equation. Notably, while the theoretical considerations are technical, the finite-dimensional representation is relatively straightforward and relies only on the quantities $M$, $\Gamma$, and $H$ from Remark~\Rref{rem:riccati}. However, in contrast to Remark~\Rref{rem:riccati}, we cannot expect that $p_i \in \clSpan\{v_1,v_2,\dots\}$, and the existence of a positive solution to the derived matrix equation is not guaranteed. 
	The quadratic matrix equation is solved by iterative linearization via Lyapunov equations, similarly to the Newton--Kleinman approach for Riccati equations~\cite{Klei68}. Let us emphasize that this procedure is not equivalent to a Newton scheme. The approximate value bilinear form is then used to derive the sum-of-squares solution using the identity given in Theorem~\Rref{thm:lowrank_optimalcontrol} \Rref{thm:lowrank_optimalcontrol:value_function}.
	For reproducibility, the implementation is available online \cite{Hoe26}.
	
	We investigate the Van der Pol oscillator, a common test example for nonlinear dynamics; see, e.g., \cite{AzmKK21,BreH23}. Similar to \cite{BreH23}, we modify the dynamics by including a friction term so that zero is the only accumulation point. To satisfy the tangent condition $f(x)^\top\nu(x)\le 0$, we add the term $-\kappa x_1^3$. We consider the domain $\Omega := [-3,3]\times[-3,3]$ and a simple quadratic cost. The resulting dynamics and cost are given as follows:
	\begin{align*}
		f_u(x_1,x_2) := \begin{pmatrix}
			x_2 - \kappa x_1^3\\
			- \mu ( x_1^2 -1) x_2 - x_1 - \eta x_2 + \gamma u(x_1,x_2)
		\end{pmatrix}, \  g(x_1,x_2 ) := x_1^2 + x_2^2 
	\end{align*}
	with $\mu = 2$, friction coefficient $\eta = 2.2$, $\kappa = 1.5\times10^{-1}$, and $\gamma = 4$; thus $b\equiv(0,\gamma)^\top$. As the weight, we again choose $w(x):=\frac{1}{\|x\|}$.
	
	With the exception of \Rref{req:exp_stab:value_function}, the assumptions of Lemma~\Rref{lemma:exponentially_stability} can be verified by simple calculations. The regularity of the value function required in \Rref{req:exp_stab:value_function} is linked to the existence of strong solutions of the HJB equation. 
	We are not aware of such existence results for infinite-horizon problems on arbitrary domains and thus do not verify the assumption analytically. Accordingly, the following experiment should be viewed as a numerical illustration conditional on
	this regularity assumption. If this unverified regularity assumption holds, the Perron--Frobenius semigroup is exponentially stable, and the results from Section~\Rref{ssec:koopman_feedback} apply.

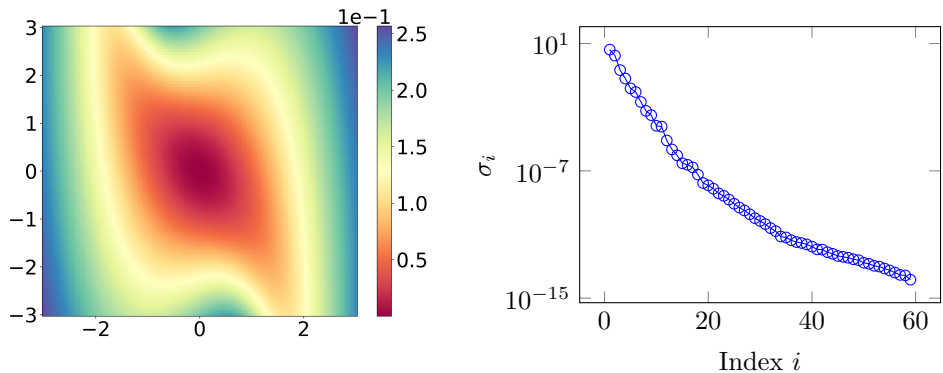
\begin{figure}[hbt]
		\begin{minipage}{.49 \textwidth}
			\includegraphics[width = \textwidth]{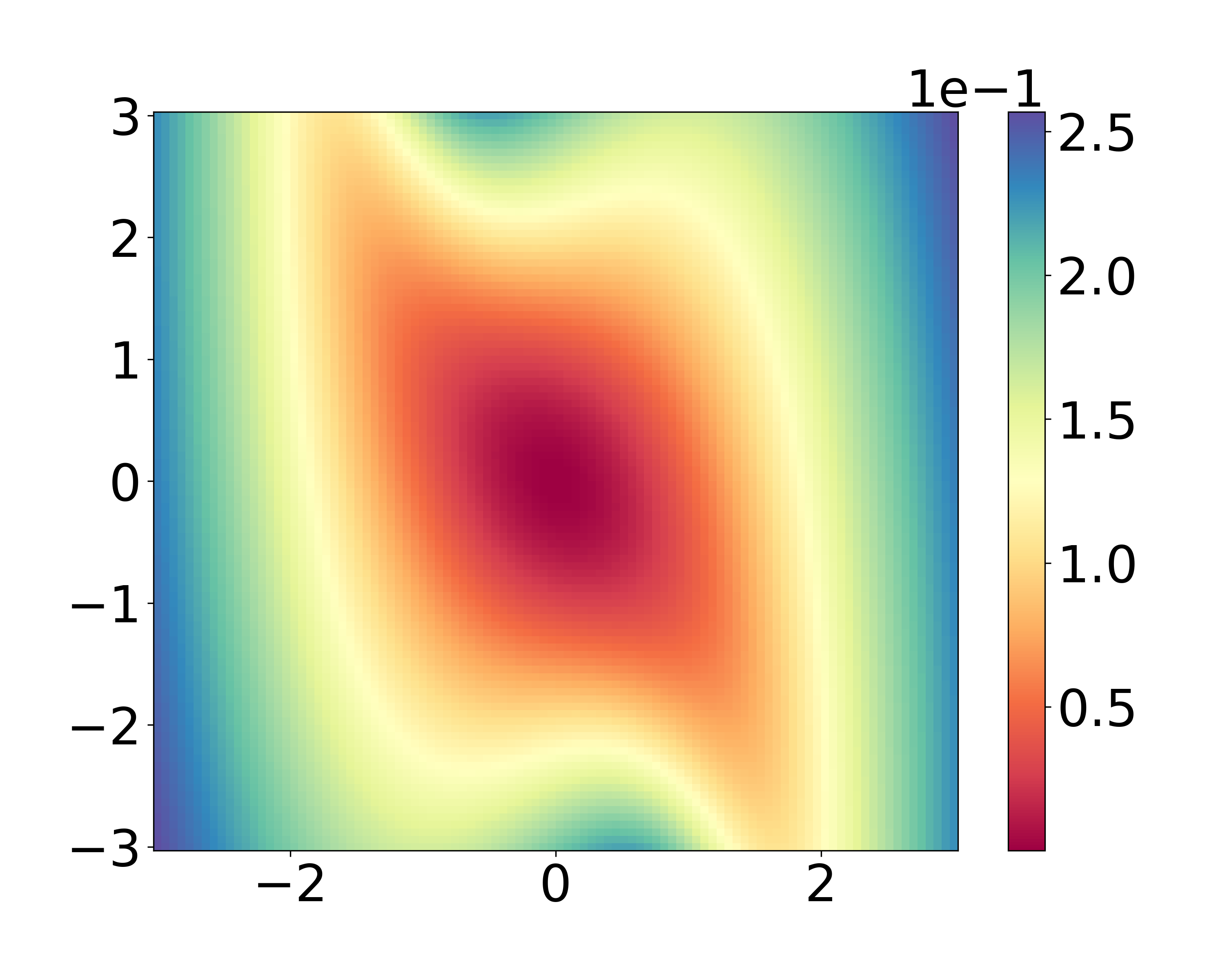}
			\vspace{.3cm}
		\end{minipage}
		\begin{minipage}{.49 \textwidth}
			\begin{tikzpicture}
				\begin{semilogyaxis}[    
					xlabel=Index $i$,    
					ylabel=$\sigma_i$,					
					legend pos=north east,
					restrict x to domain=0:59,
					width=\textwidth,
					height=.29\textheight,      
					]
					\addplot[
					mark=o,
					color=blue,
					mark options={draw=blue,fill=white}
					] table [
					col sep=space, 
					x index=0,            
					y index=1   
					] {data/pi_mu=2_eta=2.2_gamma=4.0_nx=31_px=5_ny=31_py=5/singular_values_it=9.csv};     
			\end{semilogyaxis} 
		\end{tikzpicture}
	\end{minipage} 
	\caption{Left: Computed value function $v$ for the modified Van der Pol oscillator. Right: The first 60 eigenvalues $\sigma_i$ in the eigendecomposition of the value bilinear form.}
	\label{fig:vanderpol_decay_vf}
\end{figure}

	As shown in Figure~\Rref{fig:vanderpol_decay_vf} (right), the eigenvalues of the numerical approximation to the bilinear form from Lemma~\Rref{lem:lowrank_control} decay rapidly, qualitatively consistent with Lemma~\Rref{lemma:decay}. Around the 30th basis element, however, numerical errors begin to mask this decay; by that point, the magnitudes approach the achievable precision of our ansatz space.
\appendix
\section{Technical proofs}
\label{sec:technical_proofs}
\begin{lmm}
	\label{lemma:weak_densesubset}
	Let $Z$ be a Banach space, and assume that $T^\ast(t) \in \mathcal{L}(Z^\ast)$ has the semigroup property with $T^\ast(0) = I$ and $\| T^\ast(t) \|_{\mathcal{L}(Z^\ast)} \le M$ for some $M < \infty$ and $0 \le t \le \delta$. Furthermore, let $D \subseteq Z$ be dense in $Z$. If, for every $\phi \in Z^\ast$,
	\[
	\lim_{t \rightarrow 0} \inner[Z,Z^\ast]{ z }{T^\ast(t) \phi - \phi } = 0 \qquad \text{for all} \; z \in D
	\],
	then $T^\ast$ is a weak-* continuous semigroup.
\end{lmm}
\begin{proof}
	Let $\phi \in Z^\ast$ and $z \in Z$. Then there exists some $z_n \in D$ such that $\lim_{n \rightarrow \infty} \| z - z_n \|_{Z} = 0$. In particular, for all $n \in \N$:
	\begin{align*}
		& \limsup_{t \rightarrow 0} | \inner[Z,Z^\ast]{z}{T^\ast(t) \phi - \phi } | 
		\le \limsup_{t \rightarrow 0} \big( | \inner[Z,Z^\ast]{z - z_n}{T^\ast(t) \phi - \phi }| + |\inner[Z,Z^\ast]{z_n}{T^\ast(t) \phi - \phi }|\big)\hspace{-5px}\\
		&\le  (M + 1)\| z - z_n \|_{Z} \| \phi \|_{Z^\ast} + \lim_{t \rightarrow 0 } | \inner[Z,Z^\ast]{z_n}{T^\ast(t) \phi - \phi } | 
		=   (M + 1) \| z - z_n \|_{Z} \| \phi \|_{Z^\ast}.
	\end{align*}
	Letting $n\to\infty$, we conclude that the limsup is zero and hence that $T^\ast(\cdot)$ is weak-* continuous at $0$. Since $T^\ast(\cdot)$ has the semigroup property, $T^\ast$ is weak-* continuous everywhere and is therefore a weak-* continuous semigroup.
\end{proof}
\begin{proof}[Proof of Proposition~\Rref{prop:bounds_flow}]
	\label{proof:bounds_flow}
	The first assertion follows from \cite[§12, Theorem VI]{Wal98}. For the second assertion, we consider the system
	\begin{equation}
		\label{eq:bounds_flow:ode}
	\left\{ \begin{array}{lcll}
		\tfrac{\mathrm{d}}{\mathrm{d}s} z(s;z_0) &=& - f(z(s;z_0))  \qquad & \text{for} \; s \in (0,t),\\
		z(0;z_0) &=& z_0.
	\end{array}\right. 
	\end{equation}
	For $z_0 \in \Phi^t(\Omega)$, there exists at least one $x_0$ such that $\Phi^t(x_0) = z_0$. We set 
	\[
	z(s;z_0) := \Phi^{t-s}(x_0) \in \Omega,
	\]
	which solves \eqref{eq:bounds_flow:ode}. Again, by \cite[§12, Theorem VI]{Wal98}, this solution is unique, and we define 
	\[
	\Phi^{-t} \colon \Phi^{t}(\Omega) \to  \Omega,\qquad 
	z_0 \mapsto z(t;z_0).
	\]
	By the same result, $\Phi^{-t}$ is Lipschitz continuous with constant $L$. Finally,
	$
	(\Phi^{-t} \circ \Phi^t )(z) = z.
	$
	We apply the result from \cite[Theorem 3.3]{EvaG92} and obtain 
	\[
	I =  \mathrm{D} \Phi^{-t} ( \Phi^{t}(z)) \mathrm{D} \Phi^{t}(z)
	\quad \Rightarrow  \quad   
	\mathrm{D} \Phi^{-t} ( \Phi^{t}(z)) =  ( \mathrm{D} \Phi^{t}(z) )^{-1}
	\]
	for almost every $z \in \Omega$.
	This means $\mathrm{D} \Phi^t$ is invertible almost everywhere, and we have
	\[
	( \mathrm{D} \Phi^t(\cdot) )^{-1} \in L^\infty( \Omega;\R^{d\times d})
	\quad \text{with} \; \| ( \mathrm{D} \Phi^t(\cdot) )^{-1} \|_{L^\infty} \le L .
	\]
	The norm bounds imply $L^{-d} \le |\det{\mathrm{D} \Phi^t(z)}| \le L^{d}$ almost everywhere.
\end{proof}

\begin{proof}[Proof of Proposition~\Rref{prop:explicit_Dphi}]
	\label{proof:explicit_Dphi}
	By Definition~\Rref{defi:uncontrolled_flow}, the flow satisfies
	\[
	 \Phi^t(z) = z + \int_0^t f(\Phi^s(z))\dint s \qquad \text{for a.e.}\; z \in \Omega,
	\]
	and Proposition~\Rref{prop:bounds_flow} gives $\Phi^s \in W^{1,\infty}(\Omega)$. Applying weak spatial derivatives to this identity, using the chain rule for Lipschitz continuous functions \cite[Theorem 3.3]{EvaG92} and the Fubini--Tonelli theorem \cite[Remark A6.11]{Alt12}, yields
	\begin{equation}
		\mathrm D\Phi^t(z)
		=I+\int_0^t\mathrm Df(\Phi^s(z))\mathrm D\Phi^s(z)\,\dint s.
		\label{eq:integral_Dflow}
	\end{equation}
	In particular, for almost every $z$, the map $t\mapsto\mathrm D\Phi^t(z)$ is absolutely continuous and
	\begin{equation}
		\label{eq:ode_Dflow}
		\tfrac{\mathrm d}{\mathrm dt}\mathrm D\Phi^t(z)
		=\mathrm Df(\Phi^t(z))\mathrm D\Phi^t(z)
	\end{equation}
	for almost every $t\in(0,T)$, with initial value $I$. Since $\mathrm Df(\Phi^{\cdot}(z))\in L^\infty(0,T;\R^{d\times d})$ for almost every $z$, the linear integral equation \eqref{eq:integral_Dflow} has a unique Carath\'{e}odory solution by \cite[Chapter 10, XVIII]{Wal98}.
\end{proof}
\begin{proof}[Proof of Theorem~\Rref{thm:optimal_feedback}]
	Let us define 
	\begin{align*}
	&H(z,\beta) := g(z) + \beta^2 + (f(z) + \beta b(z))^\top \nabla v(z),\\
	&g(z) := \sum_{i=1}^\infty c_i(z)^2 \qquad
	  \text{and} \qquad J_z[\alpha] := \int_0^\infty g( x(t;\alpha)) + \alpha(t)^2 \dint t.
	\end{align*}
	By nuclearity of $\{c_i\}_{i\in\N}$, we have $g \in W^{1,\infty}(\Omega)$.
	First, we show that (i) $\inf_{\beta \in \R} H(z,\beta)\ge 0$ and then (ii) $\inf_{\beta \in \R} H(z,\beta) \le 0$ to conclude equality.
	
	We start with (i). Let $\beta \in \R$ and $z \in \Omega$. Then for $\tau > 0$ small enough the solution of
	\[
	\left\{
	\begin{array}{lclcl}
		\tfrac{\mathrm{d}}{\mathrm{d}t} x_\beta(t)
		&=& f(x_\beta(t)) + \beta b(x_\beta(t))
		&\qquad& 0 < t < \tau,\\
		x_\beta(0) &=& z
	\end{array}
	\right.
	\]
	remains in $\Omega$. Let $\tilde{\alpha} \in \mathcal{A}_{x_\beta(\tau)}$ be such that
	$
	J_{x_\beta(\tau)}[\tilde{\alpha}] \le v(x_\beta(\tau)) + \varepsilon
	$
	for some $\varepsilon > 0$. We set
	\[
	\alpha(t) :=
	\left\{
	\begin{array}{lcl}
		\beta &\qquad& \text{for } 0 \le t \le \tau,\\
		\tilde{\alpha}(t-\tau) &\qquad& \text{for } t > \tau.
	\end{array}
	\right.
	\]
	Then $\alpha \in \mathcal{A}_z$ and
	\[
	v(z)
	\le J_z[\alpha]
	\le \int_0^\tau \bigl(g(x_\beta(t)) + \beta^2\bigr)\,\dint t
	+ v(x_\beta(\tau)) + \varepsilon.
	\]
	Since this holds for every $\varepsilon>0$, we conclude
	\[
	\frac{v(z)-v(x_\beta(\tau))}{\tau}
	\le
	\frac{1}{\tau}
	\int_0^\tau \bigl(g(x_\beta(t)) + \beta^2\bigr)\,\dint t.
	\]
	Sending $\tau \searrow 0$ and using $v\in C^1(\Omega)$ yields
	\[
	-(f(z)+\beta b(z))^\top \nabla v(z)
	\le g(z)+\beta^2
	\]
	and hence $H(z,\beta)\ge0$ for all $\beta\in\R$, which shows (i).
	
	For (ii), fix again $z \in \Omega$ and, for every sufficiently small $\tau>0$, choose
	$\alpha_\tau \in \mathcal{A}_z$ such that
	\[
	J_z[\alpha_\tau] \le v(z)+\tau^2.
	\]
	Since the shifted control input $t\mapsto\alpha_\tau(t+\tau)$ is admissible from
	$x(\tau;\alpha_\tau)$, we have
	\begin{align*}
		& \int_0^\tau \left(g(x(t;\alpha_\tau))+\alpha_\tau(t)^2 \right)\,\dint t
		+v(x(\tau;\alpha_\tau)) \\
		\le & \int_0^\tau \left(g(x(t;\alpha_\tau))+\alpha_\tau(t)^2 \right)\,\dint t + \int_{\tau}^\infty \left(g(x(t;\alpha_\tau)) + \alpha_\tau(t)^2\right) \dint t \le v(z) + \tau^2.
	\end{align*}
	By the chain rule, this directly leads to
	\[
	\int_0^\tau
	\left(
	g(x(t;\alpha_\tau))+\alpha_\tau(t)^2
	+\tfrac{\mathrm{d}}{\mathrm{d}t}v(x(t;\alpha_\tau))
	\right)\,\dint t
	=
	\int_0^\tau H(x(t;\alpha_\tau),\alpha_\tau(t))\,\dint t
	\le \tau^2.
	\]
	Since $g$, $f^\top\nabla v$, and $b^\top\nabla v$ are bounded, Young's inequality implies
	\[
	H(x,\beta)\ge \frac{1}{2}\beta^2-C
	\qquad
	\text{for all }x\in\Omega,\ \beta\in\R
	\]
	for some $C>0$. Hence, for $0<\tau\le1$,
	\[
	\int_0^\tau \alpha_\tau(t)^2\,\dint t
	\le C\tau.
	\]
	By Hölder's inequality, we also have
	\[
	\int_0^\tau |\alpha_\tau(t)|\,\dint t
	\le
	\tau^{1/2}
	\left(
	\int_0^\tau \alpha_\tau(t)^2\,\dint t
	\right)^{1/2}
	\le C\tau.
	\]
	Using the fact that $f$ and $b$ are bounded on $\Omega$ and
	$x(\cdot;\alpha_\tau)\subset\overline{\Omega}$, we conclude
	\begin{align*}
		\sup_{0\le t\le\tau}
		\|x(t;\alpha_\tau)-z\|
		&\le
		\int_0^\tau
		\|f(x(t;\alpha_\tau))\|\,\dint t
		+
		\int_0^\tau
		\|b(x(t;\alpha_\tau))\|
		|\alpha_\tau(t)|\,\dint t
		\le C\tau.
	\end{align*}
	Moreover, since $g$, $f^\top\nabla v$, and $b^\top\nabla v$ are Lipschitz continuous,
	\begin{align*}
		\int_0^\tau
		\left|
		H(x(t;\alpha_\tau),\alpha_\tau(t))
		-H(z,\alpha_\tau(t))
		\right|\,\dint t
		\le
		C
		\int_0^\tau
		\bigl(1+|\alpha_\tau(t)|\bigr)
		\|x(t;\alpha_\tau)-z\|\,\dint t
		\le C\tau^2.
	\end{align*}
	Therefore,
	\begin{align*}
		& \inf_{\beta\in\R} H(z,\beta)
		=
		\frac{1}{\tau}
		\int_0^\tau
		\inf_{\beta\in\R}H(z,\beta)\,\dint t
		\le
		\frac{1}{\tau}
		\int_0^\tau H(z,\alpha_\tau(t))\,\dint t\\
		\le &
		\frac{1}{\tau}
		\int_0^\tau
		H(x(t;\alpha_\tau),\alpha_\tau(t))\,\dint t
		+
		\frac{1}{\tau}
		\int_0^\tau
		\left|
		H(x(t;\alpha_\tau),\alpha_\tau(t))
		-H(z,\alpha_\tau(t))
		\right|\,\dint t
		\le C\tau.
	\end{align*}
	Sending $\tau\searrow0$ proves (ii). Hence
	\[
	0
	=
	\inf_{\beta\in\R}H(z,\beta)
	=
	\inf_{\beta\in\R}
	\left(
	g(z)+\beta^2+(f(z)+\beta b(z))^\top\nabla v(z)
	\right)
	\]
	for all $z\in\Omega$. Since the objective is strictly convex and quadratic in $\beta$, its unique minimizer is
	\[
	u_\ast(z):=-\frac{1}{2}b(z)^\top\nabla v(z)
	\]
	and the HJB equation can equivalently be written as
	\[
	g(z)+f(z)^\top\nabla v(z)
	-\frac{1}{4}\bigl(b(z)^\top\nabla v(z)\bigr)^2
	=0.
	\]
	
	Next, we show that the optimal feedback law is boundary-admissible. By the HJB equation, for every
	$z\in\Omega$ and $\beta\in\R$,
	\begin{align*}
		|\beta-u_\ast(z)|^2
		&=
		\beta^2
		+b(z)^\top\nabla v(z)\,\beta
		+\frac{1}{4}\bigl(b(z)^\top\nabla v(z)\bigr)^2\\
		&=
		g(z)+\beta^2
		+(f(z)+\beta b(z))^\top\nabla v(z)
		=
		H(z,\beta).
	\end{align*}
	Fix $z\in\Omega$ and let $\alpha_\varepsilon\in\mathcal{A}_z$ satisfy
	\[
	J_z[\alpha_\varepsilon]\le v(z)+\varepsilon.
	\]
	Then for every $T>0$, the shifted tail of $\alpha_\varepsilon$ is admissible from
	$x(T;\alpha_\varepsilon)$, and therefore
	\[
	v(x(T;\alpha_\varepsilon))
	\le
	\int_T^\infty
	\bigl(
	g(x(t;\alpha_\varepsilon))
	+\alpha_\varepsilon(t)^2
	\bigr)\,\dint t.
	\]
	Consequently,
	\begin{align*}
		& \int_0^T |\alpha_\varepsilon(t) -u_\ast(x(t;\alpha_\varepsilon))|^2\,\dint t
		=
		\int_0^T \left( g(x(t;\alpha_\varepsilon)) +\alpha_\varepsilon(t)^2 +\frac{\mathrm{d}}{\mathrm{d}t} v(x(t;\alpha_\varepsilon)) \right)\,\dint t\\
		= & \int_0^T \left( g(x(t;\alpha_\varepsilon)) +\alpha_\varepsilon(t)^2 \right)\,\dint t +v(x(T;\alpha_\varepsilon))-v(z) 
		\le J_z[\alpha_\varepsilon]-v(z)
		\le\varepsilon.
	\end{align*}
	Let $x_\ast$ denote the solution of the closed-loop equation
	\[
	\tfrac{\mathrm{d}}{\mathrm{d}t} x_\ast(t)
	=
	f(x_\ast(t))+b(x_\ast(t))u_\ast(x_\ast(t)),
	\qquad
	x_\ast(0)=z.
	\]
	Since $v\in W^{2,\infty}(\Omega)$, the map $f_{u_\ast} :=f+bu_\ast$ is Lipschitz continuous. Thus, for $0\le t\le T$,
	\begin{align*}
		\|x(t;\alpha_\varepsilon)-x_\ast(t)\|
		&\le
		C \left( \int_0^t
		\|x(s;\alpha_\varepsilon)-x_\ast(s)\|\,\dint s
		+\int_0^t
		|\alpha_\varepsilon(s)
		-u_\ast(x(s;\alpha_\varepsilon))|\,\dint s \right).
	\end{align*}
	By Gronwall's inequality and Hölder's inequality again, we conclude
	\[
	\sup_{0\le t\le T} \|x(t;\alpha_\varepsilon)-x_\ast(t)\|
	\le C_T \left( \int_0^T |\alpha_\varepsilon(t) -u_\ast(x(t;\alpha_\varepsilon))|^2\,\dint t \right)^{1/2}
	\le C_T\sqrt{\varepsilon}.
	\]
	Since $x(\cdot;\alpha_\varepsilon)\subset\overline{\Omega}$, letting $\varepsilon\searrow0$ shows that
	\[
	x_\ast(t)\in\overline{\Omega} \qquad\text{for every }t\ge0.
	\]

	Now let $z\in\partial\Omega$ and choose a sequence $z_n\in\Omega$ with $z_n\to z$. Denote by $x_\ast(\cdot;z_n)$ and
	$x_\ast(\cdot;z)$ the closed-loop solutions with the corresponding initial values.
	By continuous dependence on the initial value,
	\[
	\lim_{n \rightarrow \infty} \| x_\ast(t;z_n) - x_\ast(t;z) \| = 0
	\qquad\text{for every }t\ge0.
	\]
	Hence $x_\ast(t;z)\in\overline{\Omega}$ for every $t\ge0$. Since $\partial\Omega$ is $C^1$, let $\rho$ be a local $C^1$ defining function near $z$ such that
	\[
	\Omega \cap B_\varepsilon(z) =\{\rho<0\} \cap B_\varepsilon(z),
	\quad \text{and} \quad
	\nabla\rho(z)=|\nabla\rho(z)|\nu(z).
	\]
	Then $\rho(x_\ast(t;z))\le0$ for $t\ge0$ sufficiently small and $\rho(x_\ast(0;z))=0$. Therefore,
	\[
	0\ge\left. \tfrac{\mathrm{d}}{\mathrm{d}t} \rho(x_\ast(t;z)) \right|_{t=0}
	=
	\nabla\rho(z)^\top \left(f(z)+b(z)u_\ast(z)\right)
	\]
	and consequently
	\[
	\nu(z)^\top \left(f(z)+b(z)u_\ast(z)\right)\le0\qquad \text{for all }z\in\partial\Omega.
	\]
	Thus $u_\ast$ is boundary-admissible.
	
	It remains to prove optimality. Since $u_\ast$ is boundary-admissible, its closed-loop trajectory starting at any
	$z\in\Omega$ remains in $\overline{\Omega}$ for all $t\ge0$. Along this trajectory the HJB equation gives
	\[
	g(x_\ast(t))+u_\ast(x_\ast(t))^2+\frac{\mathrm{d}}{\mathrm{d}t}v(x_\ast(t))=0.
	\]
	Hence, for every $T>0$,
	\[
	\int_0^T\left(g(x_\ast(t))+u_\ast(x_\ast(t))^2\right)\,\dint t
	= v(z)-v(x_\ast(T)) \le v(z).
	\]
	Letting $T \rightarrow \infty$ and setting $\alpha_\ast(t):=u_\ast(x_\ast(t))$, we conclude
	$
	J_z[\alpha_\ast] \le v(z)
	$.
	Since $u_\ast$ is bounded and the closed-loop trajectory remains in $\overline{\Omega}$, we have $\alpha_\ast\in \mathrm{L}^\infty_{\mathrm{loc}}(0,\infty;\R)$ and $T_z(\alpha_\ast)=\infty$. Thus $\alpha_\ast\in\mathcal A_z$, and the definition of $v$ yields $v(z)=J_z[\alpha_\ast]$.
\end{proof}

The following result is well known under different assumptions in the literature, e.g., \cite{BreKP19, Luk69}. For our exact setup, we give a separate proof.
\begin{prpstn}
	\label{prop:stable_optimal_feedback}
	Let $(f,b)$ be a boundary-admissible controlled dynamic pair, let $\{c_i\}_{i\in\N}$ be a finite-rank nuclear sequence of observables with $c_i\equiv0$ for all $i>r$, and let $v$ be the corresponding value function from Definition~\Rref{defi:value_function}.
	Suppose that the following conditions hold:
	\begin{enumerate}[label=(\roman*)]
		\item $v \in W^{2,\infty}(\Omega)$ and $\mathrm{D}^2v$ is continuous at $0$,
		\item $c_i,f_i \in \fWpw{\infty}{w}{\Omega}$ are continuously differentiable at $0$ for $1\le i\le r$ and $1\le i\le d$, respectively.
	\end{enumerate}
	Then $v(z) = z^\top P z + o(\|z\|^2)$ as $z\to0$, where $P$ solves the algebraic Riccati equation
	\begin{equation}
		\mathrm{D}f(0)^\top P + P \mathrm{D}f(0) - P b(0)b(0)^\top P + Q = 0
		\label{eq:appendix:riccati}
	\end{equation}
	with \; $Q := \left( \nabla c_1(0), \dots, \nabla c_r(0) \right) \left( \nabla c_1(0), \dots, \nabla c_r(0) \right)^\top \in \R^{d \times d}$.
\end{prpstn}
\begin{proof}
	Let $z \in B_\varepsilon(0)$ for $\varepsilon>0$ small enough. For arbitrary $\phi\in \fWpw{\infty}{w}{\Omega}$, we have 
	\[
	\infty > \|\phi\|_{\Lpw{\infty}{w}{\Omega}} =\esssup_{z \in \Omega} |\phi(z)| w(z) \ge \limsup_{z \to 0} |\phi(z)| w(z).
	\]
	It follows that $|\phi(0)| \le \liminf_{z \to 0} \sfrac{\| \phi \|_{\Lpw{\infty}{w}{\Omega}}}{w(z)} = 0$. By differentiability, we then have
	\[
	c_i(z) = \nabla c_i(0)^\top z + o(\|z\|) \quad \text{and} \quad f(z) = \mathrm{D} f(0)  z  + o(\|z\|).
	\]
	Since $v \in \mathrm{W}^{2,\infty}(\Omega)$, $v(0) = 0$, $v \ge 0$, and $\mathrm{D}^2v$ is continuous at $0$, it follows that
	\begin{align*}
		\nabla v(0) = 0 \quad \text{and} \quad \nabla v(z) = \mathrm{D}^2 v(0) z  + o(\|z\|).
	\end{align*}
	By Theorem~\Rref{thm:optimal_feedback}, $v$ satisfies the HJB equation in the strong sense, i.e.,
	\begin{align*}
		0 &= \inf_{\alpha \in \R} \left\{ \nabla v(z)^\top \left( f(z) + b(z) \alpha \right) + g(z) + \alpha^2 \right\} \\
		&=  - \tfrac{1}{4} \left( \nabla v(z)^\top b(z) \right)^2 + \nabla v(z)^\top f(z) + g(z)\\
		&=   \tfrac{1}{2} z^\top \mathrm{D}^2 v(0) \mathrm{D}f(0) z + \tfrac{1}{2} z^\top  \mathrm{D}f(0)^\top \mathrm{D}^2 v(0) z \\
		& \quad - \tfrac{1}{4} z^\top \mathrm{D}^2 v (0) b(z) b(z)^\top\mathrm{D}^2 v (0) z + \sum_{i=1}^r z^\top \nabla c_i(0) \nabla c_i(0)^\top z + o(\| z \|^2 )
	\end{align*}
	with $g(x) := \sum_{i=1}^r c_i(x)^2$. For any fixed $y \in \R^d\setminus\{0\}$, set $z = \varepsilon y/\|y\|$. Dividing by $\varepsilon^2$, letting $\varepsilon\to0$, and multiplying by $\|y\|^2$, we conclude from the continuity of $b$ that
	\begin{align*}
		0 = & \tfrac{1}{2} y^\top \mathrm{D}^2 v(0) \mathrm{D}f(0) y + \tfrac{1}{2} y^\top  \mathrm{D}f(0)^\top \mathrm{D}^2 v(0) y \\
		& \qquad - \tfrac{1}{4} y^\top \mathrm{D}^2 v (0) b(0) b(0)^\top\mathrm{D}^2 v (0) y + \sum_{i=1}^r y^\top \nabla c_i(0) \nabla c_i(0)^\top y.
	\end{align*}
	By the polarization identity $P:= \tfrac{1}{2}\mathrm{D}^2 v(0)$ is a solution to \eqref{eq:appendix:riccati}.
\end{proof}

\begin{proof}[Proof of Proposition~\Rref{prop:nonlinear_op_eq_bd}]
	\label{pf:prop:nonlinear_op_eq_bd}
We prove the claims one by one.

\begin{enumerate}[label=(\roman*)]
	\item
	We first show that $E$ is bounded. By the definition of $E$ and since $\Psi^{-1}$ is an isometry,
	\begin{align*}
		\|E\phi\|_{\Wpw{2}{w}{\Omega}}^2
		&=
		\bigl\|\Psi^{-1}([(\phi,0,\dots,0)]_\sim)\bigr\|_{\Wpw{2}{w}{\Omega}}^2 
		=
		\|[(\phi,0,\dots,0)]_\sim\|_{\NWpw{2}{w}{\Omega}}^2 \\
		&=
		\inf_{(\tilde{\phi}_0,\dots,\tilde{\phi}_d)\sim (\phi,0,\dots,0)}
		\left(
		\|\tilde{\phi}_0\|_{\Lpw{2}{w}{\Omega}}^2
		+ \sum_{i=1}^d \|\tilde{\phi}_i\|_{\Lp{2}{\Omega}}^2
		\right) 
		\le
		\|\phi\|_{\Lpw{2}{w}{\Omega}}^2.
	\end{align*}
	Thus $E$ belongs to $\mathcal{L}(H_X;H_Y)$ and $\mathcal{L}(H_Y;H_Y)$. To prove injectivity, let $\phi \in H_X$ satisfy $E\phi = 0$. Then
	\[
	\Psi^{-1}\bigl([(\phi,0,\dots,0)]_\sim\bigr)=0,
	\]
	and since $\Psi^{-1}$ is injective, this implies
	\[
	[(\phi,0,\dots,0)]_\sim = [(0,\dots,0)]_\sim.
	\]
	The equivalence means that $\inner[H_X]{\phi}{\psi}=0$ for every $\psi\in\fWpw{\infty}{w}{\Omega}$. Since this test space is dense in $H_X$, it follows that $\phi = 0$ in $\fLpw{2}{w}{\Omega}$. Thus $E$ is injective.
	
	\item
	Since $(v_i)_{i\in\N}$ is a Riesz sequence, its synthesis operator $T$
	is bounded and injective, and its adjoint is given by
	\[
	T^\ast\colon H_Y \to \ell_2,
	\qquad
	\phi \mapsto \bigl(\inner[H_Y]{v_i}{\phi}\bigr)_{i\in\N}.
	\]
	
	Moreover, by definition of $E$,
	\[
	\inner[H_Y]{\phi}{E\psi}
	=
	\inner[H_X]{\phi}{\psi}
	=
	\inner[H_Y]{E\phi}{\psi}
	\qquad
	\text{for all } \phi,\psi \in H_Y.
	\]
	In particular, $E$ is self-adjoint on $H_Y$. Furthermore, we have
	\[
		Ma := \left( \inner[H_X]{v_i}{Ta}\right)_{i \in \N} = \left( \inner[H_Y]{v_i}{E Ta}\right)_{i \in \N} = T^\ast E T a
	\]
	and since $T$, $E$, and $T^\ast$ are bounded, so is
	$M = T^\ast E T \in \mathcal{L}(\ell_2)$. We next show that $M$ is injective. Let $a\in\ell_2$ and assume $Ma=0$. Then
	\[
	0
	=
	\inner[\ell_2]{Ma}{a}
	=
	\inner[H_Y]{ETa}{Ta}
	=
	\|Ta\|_{H_X}^2.
	\]
	Hence $Ta=0$, and since $T$ is injective, it follows that $a=0$. Thus $M$ is injective. Moreover, $M^\ast=T^\ast E^\ast T=T^\ast E T=M$, so $M$ is self-adjoint.
	\item
	Let
	\[
	B^\ast \colon H_Y \to H_X,
	\qquad
	 \phi \mapsto \tilde f^\top \nabla \phi.
	\]
	By the inequality
	\[
		\| B^\ast \phi \|^2_{H_X} = \int_\Omega (\tilde f^\top \nabla \phi)(x)^2 w^2(x) \dint x \le d^2 \max_{1 \le i \le d}  \| \tilde f_i \|^2_{\Lpw{\infty}{w}{\Omega}} \| \phi \|^2_{\Wpw{2}{w}{\Omega}}
	\]
	we conclude that the operator is bounded whenever $\tilde f_i\in X^\ast$ for $i=1,\dots,d$. For $a\in\ell_2$ and $i\in\N$, we conclude
	\begin{align*}
		(\tilde F^\ast a)_i
		&=
		\inner[H_X]{v_i}{\tilde f^\top \nabla(Ta)} 
		=
		\inner[H_Y]{v_i}{E(\tilde f^\top \nabla(Ta))} 
		=
		\bigl(T^\ast E B^\ast T a\bigr)_i.
	\end{align*}
	Hence
	$
	\tilde F^\ast = T^\ast E B^\ast T
	$.
	Since $T$, $B^\ast$, $E$, and $T^\ast$ are bounded, $\tilde F^\ast$ is bounded as a composition of bounded operators.
	The statement for $F^\ast$ and $F_k^\ast$ follows by taking $\tilde f := f$ and $\tilde f := v_k b$, respectively. Finally, if $Ta\in X^\ast$, then $(Ta)b_i\in X^\ast$ and
	\[
		\|(Ta)b_i\|_{X^\ast}\leq \|b_i\|_{\Lp{\infty}{\Omega}}\|Ta\|_{X^\ast},
	\]
	so the same argument with $\tilde f:=(Ta)b$ proves that $F[a]^\ast\in\mathcal L(\ell_2)$.
	\item
	By definition,
	$
	Ha = \bigl(\inner[H_X]{c_i}{Ta}\bigr)_{i\in\N}
	$ for $a\in\ell_2$. Therefore,
	\begin{align*}
		\|Ha\|_{\ell_2}^2
		&=
		\sum_{i=1}^\infty |\inner[H_X]{c_i}{Ta}|^2 
		= \sum_{i=1}^\infty |\inner[H_Y]{c_i}{ETa}|^2\\
		\le & C
		\left( \sum_{i=1}^\infty \|c_i\|_{Y^\ast}^2 \right) \|E\|^2_{\mathcal{L}(H_Y;H_Y)} \| T \|^2_{\mathcal{L}(\ell_2;H_Y)} \| a\|^2_{\ell_2} 
	\end{align*}
	Here we used the continuous embedding $Y^\ast\subset H_Y$ and the fact that $\sum_{i=1}^\infty \|c_i\|_{Y^\ast}^2  < \infty$, since $\{c_i\}_{i \in \N}$ is a nuclear sequence of observables. We conclude that $H\in\mathcal{L}(\ell_2)$.
\end{enumerate}
\end{proof}
\paragraph{\textbf{Author contributions.}}
Bernhard Höveler was responsible for the conceptualization, methodology, investigation, implementation, validation, visualization, and preparation of the original manuscript. Tobias Breiten contributed to validation, funding acquisition, and supervision. Both authors contributed to reviewing and editing the manuscript and approved the final version.
\paragraph{\textbf{Generative AI statement:}}
OpenAI's GPT-5.6 Sol was used for language editing and proofreading. All mathematical results presented in this paper are the original work of the authors. All AI-assisted suggestions were independently and carefully verified by the authors, who take full responsibility for the content of the manuscript. 
\paragraph{\textbf{Data availability statement:}}
The exact Python implementation used in Section~\ref{sec:numerics} and to generate Figure~\Rref{fig:vanderpol_decay_vf} is available at Zenodo \cite{Hoe26}.
\bibliographystyle{plain}
\bibliography{references}
\end{document}